\documentclass[10pt]{article}
\usepackage{lmodern}
\usepackage{amsmath}
\usepackage[T1]{fontenc}
\usepackage[utf8]{inputenc}
\usepackage{authblk}
\usepackage{amsfonts}
\usepackage{graphicx}
\usepackage{rotating}
\usepackage{amssymb}
\usepackage[english]{babel}
\usepackage{xcolor}
\usepackage{amsthm}
\usepackage{graphicx}
\usepackage{mathrsfs}
\usepackage{makecell}
\usepackage{microtype}
\usepackage{mathscinet}
\usepackage{array}
\usepackage{multirow}
\usepackage{enumerate}
\usepackage[cal=boondoxo,bb=ams]{mathalfa}
\usepackage{hyperref}
\hypersetup{hidelinks}

\newtheorem{defn}{Definition}[section]
\newtheorem{theorem}{Theorem}[section]

\newtheorem{lemma}{Lemma}[section]
\newtheorem{coro}{Corollary}[section]
\newtheorem{remark}{Remark}[section]
\newtheorem{exam}{Example}[section]

\newcommand{\ml}{\mathcal}
\newcommand{\mb}{\mathbb}
\DeclareMathOperator{\lin}{lin}
\DeclareMathOperator{\nlin}{nlin}

\def\XXint#1#2#3{{\setbox0=\hbox{$#1{#2#3}{\int}$ }
		\vcenter{\hbox{$#2#3$ }}\kern-.6\wd0}}

\title{On the Cauchy problem for semilinear thermoelastic plate systems\\ in the $L^q$ framework}
\author[1]{Halit Sevki Aslan\thanks{Halit Sevki Aslan (halitsevkiaslan@gmail.com)}}
\affil[1]{Department of Computer Science and Mathematics, University of S\~ao Paulo (USP), \authorcr  14040-901 Ribeir\~ao Preto, SP, Brazil}

\author[2]{Wenhui Chen\thanks{Wenhui Chen (wenhui.chen.math@gmail.com)}}
\affil[2]{School of Mathematics and Information Science, Guangzhou University,\authorcr 510006 Guangzhou, China}

\date{}

\begin{document}
		\maketitle

		\begin{abstract}
			\medskip
		We mainly consider semilinear thermoelastic plate systems with general power nonlinearities in the whole space $\mathbb{R}^n$. By applying the Fourier analysis, some sharp $(L^q\cap L^m)-L^q$ estimates of solutions (with any $1\leqslant m\leqslant q\leqslant +\infty$) to the classical thermoelastic plate system are derived, which cover all known results in $\mathbb{R}^n$. Then, we investigate global in time existence of small data $L^q$ solutions (with any $q\in[1,+\infty]$) and blow-up of weak solutions for the semilinear thermoelastic plate systems under suitable conditions on the power exponents, which justify critical exponents for several classes of nonlinearities.
			\\
			
			\noindent\textbf{Keywords:} thermoelastic plate system, semilinear Cauchy problem, power nonlinearity, global in time existence, blow-up, critical exponent\\
			
			\noindent\textbf{AMS Classification (2020)} 35G55, 35G40, 35B40, 35Q79, 35A01
		\end{abstract}
\fontsize{12}{15}
\selectfont

\section{Introduction}\label{Section_Introduction}\setcounter{equation}{0}
\hspace{5mm}In this paper, we mainly consider the following Cauchy problem for the semilinear thermoelastic plate systems with general nonlinearities:
\begin{align}\label{General-Semilinear-TEP}
\begin{cases}
u_{tt}^N+\Delta^2u^N+\Delta\theta^N=f_1^N(u^N,|\nabla|^{\sigma}u^N,u_t^N,\theta^N;p_1),&x\in\mb{R}^n,\ t>0,\\
\theta^N_t-\Delta\theta^N-\Delta u_t^N=f_2^N(u^N,|\nabla|^{\sigma}u^N,u_t^N,\theta^N;p_2),&x\in\mb{R}^n,\ t>0,\\
(u^N,u_t^N,\theta^N)(0,x)=(u^N_0,u^N_1,\theta^N_0)(x),&x\in\mb{R}^n,
\end{cases}
\end{align}
with the power nonlinearities $f_j^N=f_j^N(u^N,|\nabla|^{\sigma}u^N,u_t^N,\theta^N;p_j)$ such that
\begin{align}\label{nonlinearity}
f_j^N:=c_{1,j}|u^N|^{p_j}+c_{2,j}|\,|\nabla|^{\sigma}u^N|^{p_j}+c_{3,j}|u_t^N|^{p_j}+c_{4,j}|\theta^N|^{p_j}\ \ \mbox{with}\ \ c_{d,j}\geqslant0,
\end{align}
 where $\sigma\in\{1,2\}$ and $p_j>1$.
 The operator $|\nabla|^{\sigma}$ on the right-hand sides with $\sigma>0$ is defined via
 \begin{align*}
 |\nabla|^{\sigma}g:=\ml{F}^{-1}(|\xi|^{\sigma}\widehat{g}\,)\ \ \mbox{for all}\ \ g\in H^{-\infty},
 \end{align*}
so that $|\nabla|^2g=-\Delta g$. The functions with the superscript ``$N$'' mean with emphasis the components for the nonlinear models.
  In particular, taking $c_{d,j}\equiv0$ for all $d\in\{1,\dots,4\}$ and $j\in\{1,2\}$, the Cauchy problem \eqref{General-Semilinear-TEP} reduces to its corresponding linearized model, i.e. the classical thermoelastic plate system with the Fourier law of heat conduction, namely,
\begin{align}\label{Linear-TEP}
	\begin{cases}
		u_{tt}+\Delta^2u+\Delta\theta=0,&x\in\mb{R}^n,\ t>0,\\
		\theta_t-\Delta\theta-\Delta u_t=0,&x\in\mb{R}^n,\ t>0,\\
		(u,u_t,\theta)(0,x)=(u_0,u_1,\theta_0)(x),&x\in\mb{R}^n,
	\end{cases}
\end{align}
where $u=u(t,x)$ and $\theta=\theta(t,x)$ denote, respectively, the vertical displacement of plate and the temperature difference (relative to some reference temperature). Our purpose in the present manuscript is to derive some qualitative properties of solutions to the Cauchy problems \eqref{General-Semilinear-TEP} and \eqref{Linear-TEP} in the $L^q$ framework with any $q\in[1,+\infty]$, including the global in time existence and sharp estimates of solutions. As a corollary, the global in time existence result is also generalized to strongly coupled nonlinearities of the product-type. Then, by investigating blow-up of weak solutions in finite time under suitable conditions for $p_1,p_2$ in some classes of power nonlinearities \eqref{nonlinearity}, some critical exponents for the nonlinear Cauchy problem \eqref{General-Semilinear-TEP} are shown.

\subsection{Background for thermoelastic plate systems}
\hspace{5mm}It is well-known that the classical thermoelastic plate system (i.e. the plate or beam equation coupled with the Fourier law of heat conduction) describes a homogeneous, elastic and thermally isotropic plate subjecting to a temperature distribution, which has been deeply studied in the last forty years (see \cite{Lagnese-Lions=1988,Lagnese=1989,Kim=1992,Munoz-Racke=1995,Liu-Zheng=1997} and numerous references therein). Particularly, concerning its corresponding Cauchy problem, i.e. the linear coupled system  \eqref{Linear-TEP}, the pioneering works \cite{Racke=1992,Munoz-Racke=1995} derived some $L^2$ estimates for the classical energy term
\begin{align*}
\mb{E}(t,x):=\big(u_t(t,x),\Delta u(t,x),\theta(t,x)\big)^{\mathrm{T}}
\end{align*}
via multiplicative techniques in the Fourier space. Later, by using the Riesz-Thorin interpolation argument between $L^2\mapsto L^2$ and $L^1\mapsto L^{\infty}$, when $q\in[2,+\infty]$ a $L^q$ decay estimate for the energy term $\mb{E}(t,\cdot)$ with the $L^{q'}$ initial data $\mb{E}(0,\cdot)$ were obtained in \cite{Denk-Racke=2006}. According to the Mikhlin-H\"ormander multiplier theorem, some resolvent estimates in the $L^q$ framework with $q\in(1,+\infty)$ were found in \cite{Naito-Shibata=2009,Denk-Racke-Shibata=2009}. By the first order expansion of $L^1$ functions (see \cite[Lemma 3.1]{Ikehata=2004} originally and \cite[Lemma 5.1]{Ikehata-Michihisa=2019} for the higher order version), the author of \cite{Said-Houari=2013} got a faster  decay $L^2$ estimate of the energy term $\mb{E}(t,\cdot)$ with the $(1+|x|)^{\gamma}$-weighted $L^1$ initial data satisfying $\int_{\mb{R}^n}\mb{E}(0,x)\,\mathrm{d}x=0$. To investigate sharp $L^2$ estimates of $\mb{E}(t,\cdot)$, by employing energy methods in the Fourier space as well as explicit characteristic roots, the paper \cite{Racke-Ueda=2016} derived a sharp $(L^2\cap L^1)-L^2$ decay estimate for the energy term. Recently, the authors of \cite{Chen-Ikehata=2023} discovered optimal growth estimates (when $n=1,2,3$); optimal neither growth nor decay estimate (when $n=4$); optimal decay estimates (when $n\geqslant 5$) for the vertical displacement $u(t,\cdot)$ in the $L^2$ norm as large time, which implies the critical dimension $n=4$ for distinguishing the decisive role (for large time behavior) between the pure plate model and the Fourier law of heat conduction thanks to the diffusion plate asymptotic profile. Furthermore, singular limits with respect to the small thermal parameter in the Fourier law were deduced. Nevertheless, sharp estimates for the solutions $u(t,\cdot)$ and $\theta(t,\cdot)$ themselves, even for the standard energy term $\mb{E}(t,\cdot)$, in the general $L^q$ norm with  $q\in[1,+\infty]$ are still open, especially at the endpoints $q\in\{1,+\infty\}$, due to some technical difficulties (see the introduction in Section \ref{Section-Linear} for detailed explanations). We in Theorem \ref{Thm-Linear-u} are going to completely answer this question by the means of the Fourier analysis. 

In recent years, several corresponding initial-boundary value problems (with general nonlinearities) to the semilinear thermoelastic plate systems \eqref{General-Semilinear-TEP}, namely,
\begin{align}\label{Boundary-Value-Problem}
\begin{cases}
u_{tt}^N+\Delta^2u^N+\Delta\theta^N=f_1^N(u^N,|\nabla|^{\sigma}u^N,u_t^N,\theta^N),&x\in\Omega,\ \,\,\, t>0,\\
\theta^N_t-\Delta\theta^N-\Delta u_t^N=f_2^N(u^N,|\nabla|^{\sigma}u^N,u_t^N,\theta^N),&x\in\Omega,\ \,\,\, t>0,\\
u^N=\Delta u^N=\theta^N=0,&x\in\partial\Omega,\ \! \, t>0,\\
(u^N,u_t^N,\theta^N)(0,x)=(u^N_0,u^N_1,\theta^N_0)(x),&x\in\Omega,
\end{cases}
\end{align}
where $\Omega\subset\mb{R}^n$ is a bounded domain with smooth enough boundary $\partial\Omega$, have caught a lot of attentions. We refer the interested reader to: \cite{Leiva-Sivoli=2003} for the existence, stability and smoothness of a bounded solution to the model \eqref{Boundary-Value-Problem} with $f_{j}^N(u^N,\theta^N)$ carrying $j\in\{1,2\}$; \cite{Naso=2005} for the controllability to the model \eqref{Boundary-Value-Problem} with $f_1^N(u^N_t)$ and $f_2^N(\theta^N)$; \cite{Baroun-Boulite-Diag-Maniar=2009} for the  existence of almost periodic solutions to the model \eqref{Boundary-Value-Problem} with $f_{j}^N(\nabla u^N,\nabla\theta^N)$ carrying $j\in\{1,2\}$; \cite{Fatori-Jor-Ma-Yang=2015,Bezerra-Carbone-Nasci-Schiabel=2018,Bezerra-Carbone-Nascimento-Schiabel=2019} for the pullback attractors or the global exponential attractors to the model \eqref{Boundary-Value-Problem} with $f_1^N(u^N)$ and $f_2^N=0$; \cite{Singh-Vijayakumar-Shukla-Chuahan=2022} for the asymptotic stability to the model \eqref{Boundary-Value-Problem} with $f_1^N(u^N)$ and $f_2^N(\theta^N)$. For another, concerning related nonlinear Cauchy problems for the thermoelastic plate system \eqref{Linear-TEP}, one may see \cite{Racke-Ueda=2017,Banquet-Doria-Villamizar-Roa=2024} and references given therein. To the best of our knowledge, there is no any result for  semilinear thermoelastic plate systems in the whole space $\mb{R}^n$, particularly, for the general semilinear systems \eqref{General-Semilinear-TEP}. The difficulties of studying our nonlinear problem \eqref{General-Semilinear-TEP} in the general $L^q$ framework are
\begin{itemize}
	\item the invalidity of Poincar\'e inequality in $\mb{R}^n$ (comparing with boundary value problems);
	\item the general index $q\in[1,+\infty]$ (comparing with the $L^2$ or energy frameworks);
	\item the mixed nonlinearities consisting of $u^N$ and $\theta^N$, even with $u^N_t$, $|\nabla|u^N$, $\Delta u^N$;
	\item the blow-up mechanism from the nonlinearities.
\end{itemize}

\subsection{Critical exponents for semilinear structurally damped plate equations}
\hspace{5mm} Let us state our observation on critical exponents for the semilinear thermoelastic plate systems \eqref{General-Semilinear-TEP} with the nonlinearities \eqref{nonlinearity} carrying $c_{d,2}=0$ for $d\in\{1,\dots,4\}$.
Neglecting the time variation of temperature difference unknown $\theta^N$ formally and taking $f_2^N=0$ in our Cauchy problem \eqref{General-Semilinear-TEP}, one may find $\Delta \theta^N=-\Delta u_t^N$, which reduces  \eqref{General-Semilinear-TEP} to the following (scalar) semilinear structurally damped plate  equations:
\begin{align}\label{Semilinear-Sigma}
\begin{cases}
u_{tt}^N+\Delta^2u^N-\Delta u_t^N=f_1^N(u^N,u_t^N;p),&x\in\mb{R}^n,\ t>0,\\
(u^N,u_t^N)(0,x)=(u_0^N,u_1^N)(x),&x\in\mb{R}^n,
\end{cases}
\end{align}
in which we took $c_{2,1}=c_{4,1}=0$ such that $f_1^N(u^N,u_t^N;p)=c_{1,1}|u^N|^{p}+c_{3,1}|u^N_t|^{p}$.

The critical exponents $p=p_{\mathrm{crit}}(n)$ for \eqref{Semilinear-Sigma} with the $L^1$ data are well-understood as
\begin{align}\label{Critical-sigma}
p_{\mathrm{crit}}(n)=\begin{cases}
\displaystyle{\bar{p}_c(n):=1+\frac{4}{(n-2)_+}}&\mbox{if}\ \ c_{1,1}>0,\ c_{3,1}=0,\\[1em]
\displaystyle{p_c(n):=1+\frac{2}{n}}&\mbox{if}\ \ c_{1,1}=0,\ c_{3,1}>0,
\end{cases}
\end{align}
whose rigorous justifications are given by \cite{Pham-Kainae-Reissig=2015,D'Abbicco-Ebert=2017,Ebert-Reissig=2018,Dao=2022}.
By the critical exponent $p_{\mathrm{crit}}(n)$ here we mean a threshold value in the range for the power of nonlinear term $|\partial_t^ku^N|^{p}$ such that for $1<p\leqslant p_{\mathrm{crit}}(n)$ local in time solutions blow up in finite time under suitable sign assumptions for the Cauchy data, while for $p>p_{\mathrm{crit}}(n)$ solutions are globally in time defined (in a suitable function space) provided that the Cauchy data are sufficiently small.

 Turning back to our problem \eqref{General-Semilinear-TEP}, the coupling structures between the vertical displacement $u^N$ and the temperature difference $\theta^N$, appearing in not only the linear part but also the nonlinear part on the right-hand sides, generate additional difficulties. Naturally, one may further wonder critical exponents for the semilinear coupled systems \eqref{General-Semilinear-TEP} with the general power nonlinearities \eqref{nonlinearity}.

\subsection{Main purposes of this manuscript}
\hspace{5mm}
Our first main goal is to study sharp $L^q$ estimates of solutions to the classical thermoelastic plate system \eqref{Linear-TEP} with the $L^q\cap L^m$ initial data and $1\leqslant m\leqslant q\leqslant +\infty$, which cover all known estimates in \cite{Munoz-Racke=1995,Denk-Racke=2006,Racke-Ueda=2016,Chen-Ikehata=2023} for $\mb{R}^n$ (see comparison in Section \ref{Section-Linear}). Precisely, by suitable representations of solutions and the Fourier analysis, in Section \ref{Section-Linear} we overcome strong singularities by oscillations to derive sharp $L^r$ estimates (with any $r\in[1,+\infty]$) for the kernels implying our desired $L^q$ estimates of solutions. As byproducts, the $L^q$ well-posedness and the ultra-analytic smoothing effect are derived.

Then, by constructing a suitable time-weighted Sobolev space and applying the Banach fixed point argument, in Section \ref{Section-GESDS} we demonstrate the global in time existence result as well as the sharp $L^q$ estimates of solutions to the semilinear thermoelastic plate systems \eqref{General-Semilinear-TEP} under some conditions on the powers $p_1,p_2$ in \eqref{nonlinearity}, where the global in time solutions belong to the $L^q$ space with any $q\in[1,+\infty]$. As a corollary, we claim the global in time existence result for the model with stronger coupled nonlinearities of the product-type \eqref{strong-nonlinearity} than the additivity-type \eqref{nonlinearity}. Finally, the critical exponents for several classes of nonlinearities (see Examples \ref{Exam-01} and \ref{Exam-02}) are found by proving the blow-up of weak solutions in Section \ref{Section-Blow-up} under sign conditions for the initial data via test function methods.

\paragraph{\large Notation} The constants $c$ and $C$ may be changed from line to line. We write $g\lesssim h$ if there exists a positive constant $C$ such that $g\leqslant Ch$. We write $(a)_+:=\max\{a,0\}$ and, as usual, $1/(a)_+:=+\infty$ if $a\leqslant0$. The H\"older's conjugate of $q\in[1,+\infty]$ is taken by $q'=q/(q-1)$. We denote by $H^s_q$ and $\dot{H}^s_q$, respectively, the Bessel potential space and the Riesz potential space with $s\in\mb{R}$ regularity based on the Lebesgue space $L^q$ with $q\in[1,+\infty]$. The differential operator $|D|^{s}$ has its symbol $|\xi|^s$ with any $s\in\mb{R}$.

\section{Main results}\setcounter{equation}{0}\label{Section-Main-Result}
\subsection{Global in time existence of solutions}
\hspace{5mm}Let us first introduce the initial data space for the sake of convenient
\begin{align*}
\ml{D}:=(H^2_1\cap H^2_{\infty})\times (L^1\cap L^{\infty})\times (L^1\cap L^{\infty}).
\end{align*}
Concerning the nonlinear Cauchy problem \eqref{General-Semilinear-TEP} carrying the power nonlinearities \eqref{nonlinearity} of the additivity-type, we state the next global in time existence result.
\begin{theorem}\label{Thm-GESDS}
Let us consider the nonlinearities \eqref{nonlinearity} with $\sigma\in\{1,2\}$ and their powers
\begin{align}\label{GESDS-Condition}
	p_j>\begin{cases}
		\displaystyle{1+\frac{4}{n-2}}&\mbox{if}\ \ c_{1,j}\neq0,\ \ \mbox{and}\ \ n\geqslant 3,\\[1em]
		\displaystyle{1+\frac{4-\sigma}{n+\sigma-2}}&\mbox{if}\ \ c_{1,j}=0\ \ \mbox{but}\ \ c_{2,j}\neq0,\ \ \mbox{and}\ \ n\geqslant 3-\sigma\\[1em]
		\displaystyle{1+\frac{2}{n}}&\mbox{if}\ \ c_{1,j}=c_{2,j}=0\ \ \mbox{but}\ \ c_{3,j}\neq0\ \ \mbox{or} \ \ c_{4,j}\neq0,\\[1em]
		1&\mbox{if}\ \ c_{1,j}=c_{2,j}=c_{3,j}=c_{4,j}=0,
	\end{cases}
\end{align}
for $j\in\{1,2\}$. Suppose that the initial data $(u_0^N,u_1^N,\theta_0^N)\in\ml{D}$. Then, there exists $\epsilon>0$ such that for all $\|(u_0^N,u_1^N,\theta_0^N)\|_{\ml{D}}\leqslant \epsilon$,
there is a pair of uniquely determined Sobolev solutions
\begin{align}\label{Solution-space-Lq}
	(u^N,\theta^N)\in\big(\ml{C}([0,+\infty),H^2_q)\cap\ml{C}^1([0,+\infty),L^q)\big)\times \ml{C}([0,+\infty),L^q)
\end{align}
for any $q\in[1,+\infty]$, to the semilinear thermoelastic plate systems \eqref{General-Semilinear-TEP}. Furthermore, the following sharp estimates of solutions in the $L^q$ framework with $q\in[1,+\infty]$:
\begin{align}
	\|\,|D|^s\partial_t^k u^N(t,\cdot)\|_{L^q}&\lesssim(1+t)^{-\frac{n}{2}\left(1-\frac{1}{q}\right)-\frac{s}{2}-k+1}\|(u_0^N,u_1^N,\theta_0^N)\|_{\ml{D}},\label{Solution-uN-Est}\\
	\| \theta^N(t,\cdot)\|_{L^q}&\lesssim(1+t)^{-\frac{n}{2}\left(1-\frac{1}{q}\right)}\|(u_0^N,u_1^N,\theta_0^N)\|_{\ml{D}},\label{Solution-thetaN-Est} 
\end{align}
hold for any $t\geqslant 0$ with $s\in\{0,1,2\}$ if $k=0$, and $s=0$ if $k=1$.
\end{theorem}
\begin{remark}
As we will show later, the time-dependent functions in the $L^q$ estimates of Theorem \ref{Thm-GESDS} exactly coincide with those of the corresponding linearized thermoelastic system \eqref{Linear-TEP}, i.e. Theorem \ref{Thm-Linear-u} with $m=1$, which verifies the effect of no loss of estimates.
\end{remark}
\begin{remark}\label{Rem-1-infty}
Concerning $q\in(1,+\infty)$, thanks to the well-known fractional Gagliardo-Nirenberg inequality in the Riesz potential space based on $L^q$ (see \cite{Hajaiej-Molinet-Ozawa-Wang-2011} and references therein) one finds
\begin{align*}
\|\,|D|^{\kappa}u^N(t,\cdot)\|_{L^q}\lesssim \|u^N(t,\cdot)\|_{L^q}^{1-\frac{\kappa}{2}}\|\,|D|^2u^N(t,\cdot)\|_{L^{q}}^{\frac{\kappa}{2}}
\end{align*}
for any $\kappa\in(0,2)$. Then, the estimate \eqref{Solution-uN-Est} can be improved for any $s\in[0,2]$. For another, with the same philosophy, Theorem \ref{Thm-GESDS} can be generalized to the nonlinearities \eqref{nonlinearity} with $\sigma\in(0,2]$ if we consider $L^q$ solutions with $q\in[1,+\infty)$.
\end{remark}

Indeed, thanks to the sharp $L^q$ estimates in Theorem \ref{Thm-Linear-u} for the linear Cauchy problem \eqref{Linear-TEP}, we are able to study the semilinear thermoelastic plate systems \eqref{General-Semilinear-TEP} with the following strongly coupled nonlinearities:
\begin{align}\label{strong-nonlinearity}
f_1^N=|u^N|^{p_1}\,|\,|\nabla|^{\sigma}u^N|^{p_2}\,|u^N_t|^{p_3}\,|\theta^N|^{p_4}\ \ \mbox{and} \  \ f_2^N=|u^N|^{\bar{p}_1}\,|\,|\nabla|^{\sigma}u^N|^{\bar{p}_2}\,|u^N_t|^{\bar{p}_3}\,|\theta^N|^{\bar{p}_4}
\end{align}
with $p_{\ell},\bar{p}_{\ell}\in\{0\}\cup(1,+\infty)$ for $\ell\in\{1,\dots,4\}$. The last strongly coupled nonlinearities are generalizations of the weakly coupled nonlinearities, e.g. $f_1^N=|u_t^N|^{p_3}$ and $f_2^N=|\theta^N|^{\bar{p}_4}$ if other powers are equal to zero. This philosophy is also valid for studying other complicated nonlinearities including $u^N,\nabla u^N,u_t^N,\Delta u^N,\theta^N$, which is beyond the scope of
our paper
\begin{coro}\label{Coro-GESDS=2}
	Let us consider the nonlinearities \eqref{strong-nonlinearity} with $\sigma\in\{1,2\}$ and their powers
	\begin{align*}
	p_1+\left(1-\frac{\sigma}{2}\right)p_2-\frac{n}{2}\sum\limits_{\ell=1,\dots,4}p_{\ell}&<-1-\frac{n}{2},\\
	\bar{p}_1+\left(1-\frac{\sigma}{2}\right)\bar{p}_2-\frac{n}{2}\sum\limits_{\ell=1,\dots,4}\bar{p}_{\ell}&<-1-\frac{n}{2}.
	\end{align*}
	Suppose that the initial data $(u_0^N,u_1^N,\theta_0^N)\in\ml{D}$. Then, there exists $\epsilon>0$ such that for all $\|(u_0^N,u_1^N,\theta_0^N)\|_{\ml{D}}\leqslant \epsilon$,
	there is a pair of uniquely determined Sobolev solutions \eqref{Solution-space-Lq} for any $q\in[1,+\infty]$, to the semilinear thermoelastic plate systems \eqref{General-Semilinear-TEP}. Furthermore, the global in time solutions satisfy \eqref{Solution-uN-Est} and \eqref{Solution-thetaN-Est} for any $t\geqslant 0$ with $s\in\{0,1,2\}$ if $k=0$, and $s=0$ if $k=1$.
\end{coro}

\subsection{Blow-up of solutions in finite time}
\hspace{5mm}Let us turn to blow-up phenomena for the semilinear thermoelastic plate systems  under some contrary conditions to \eqref{GESDS-Condition} in the power nonlinearities \eqref{nonlinearity}, where we concentrate on the case $c_{4,1}+c_{4,2}>0$ (i.e. the coefficients of $|\theta^N|^{p_j}$ in the nonlinearities $f_j^N$ are non-trivial simultaneously).  To clarify some influence of nonlinearities in regard to $u^N$ we next separate our discussion into two parts: $c_{1,1}+c_{1,2}=0$ and $c_{1,1}+c_{1,2}>0$.

For one thing, concerning $c_{1,1}+c_{1,2}=0$, i.e. $c_{1,1}=c_{1,2}=0$ leading to the vanishing coefficients of $|u^N|^{p_j}$, we first introduce a suitable notation of weak solution to the semilinear thermoelastic plate systems \eqref{General-Semilinear-TEP}.
\begin{defn}\label{Defn-weak-01}
	A pair of functions $(u^N,\theta^N)$ is said to be a global in time weak solution to the semilinear Cauchy problem \eqref{General-Semilinear-TEP} with the nonlinearities \eqref{nonlinearity} carrying $c_{1,1}+c_{1,2}=0$ if
	\begin{align*}
		\begin{cases}
	|\nabla |^{\sigma}u^N\in L^{p_j}_{\mathrm{loc}}([0,T)\times \mb{R}^n)&\mbox{if}\ \ c_{2,j}>0,\\
	u^N_t\in L^{p_j}_{\mathrm{loc}}([0,T)\times \mb{R}^n)&\mbox{if}\ \ c_{3,j}>0,\\
	\theta^N\in L^{p_j}_{\mathrm{loc}}([0,T)\times \mb{R}^n)&\mbox{if}\ \ c_{4,j}>0,		
		\end{cases}
	\end{align*}
and	 the following integral equalities:
	\begin{align*}
		&\int_0^{+\infty}\int_{\mb{R}^n}f_1^N(t,x)\Psi_1(t,x)\,\mathrm{d}x\,\mathrm{d}t+\int_{\mb{R}^n}\Big(-u^N_0(x)\int_0^{+\infty}\Delta^2\Psi_1(\tau,x)\,\mathrm{d}\tau+u^N_1(x)\Psi_1(0,x)\Big)\,\mathrm{d}x\\
		&=\int_0^{+\infty}\int_{\mb{R}^n}u^N_t(t,x)\Big(-\partial_t\Psi_1(t,x)+\int_t^{+\infty}\Delta^2\Psi_1(\tau,x)\,\mathrm{d}\tau\Big)\,\mathrm{d}x\,\mathrm{d}t+\int_0^{+\infty}\int_{\mb{R}^n}\theta^N(t,x)\Delta \Psi_1(t,x)\,\mathrm{d}x\,\mathrm{d}t
	\end{align*}
as well as
	\begin{align*}
		&\int_0^{+\infty}\int_{\mb{R}^n}f_2^N(t,x)\Psi_2(t,x)\,\mathrm{d}x\,\mathrm{d}t+\int_{\mb{R}^n}\theta_0^N(x)\Psi_2(0,x)\,\mathrm{d}x\\
		&=\int_0^{+\infty}\int_{\mb{R}^n}\theta^N(t,x)\big(-\partial_t\Psi_2(t,x)-\Delta \Psi_2(t,x)\big)\,\mathrm{d}x\,\mathrm{d}t-\int_0^{+\infty}\int_{\mb{R}^n}u^N_t(t,x)\Delta \Psi_2(t,x)\,\mathrm{d}x\,\mathrm{d}t
	\end{align*}
	hold for any test functions $\Psi_1,\Psi_2\in\ml{C}_0^{\infty}([0,+\infty)\times\mb{R}^n)$.
\end{defn}

\begin{theorem}\label{Thm-Blow-up-01}
	Let us consider the nonlinearities \eqref{nonlinearity} with $\sigma=2$ satisfying
	\begin{itemize}
		\item $c_{d,j}\geqslant0$ for any $d\in\{1,\dots,4\}$ and $j\in\{1,2\}$;
		\item $c_{1,1}+c_{1,2}=0$,  $c_{3,1}+c_{3,2}>0$, $c_{4,1}+c_{4,2}>0$;
	\end{itemize}
and their powers 
	\begin{align}\label{Blow-up-Condition-01}
		1<p_j\leqslant 1+\frac{2}{n} \ \ \mbox{if}\ \ c_{3,j}>0\ \ \mbox{or}\ \ c_{4,j}>0
	\end{align}
	for $j\in\{1,2\}$. Suppose that the initial data $u_0^N=0$ and $u_1^N+\theta_0^N\in L^1$ verifying
		\begin{align*}
	\int_{\mb{R}^n}\big(u_1^N(x)+\theta_0^N(x)\big)\,\mathrm{d}x>0.
	\end{align*}
Then, every local in time weak solution in the sense of Definition \ref{Defn-weak-01} to the general semilinear thermoelastic plate systems \eqref{General-Semilinear-TEP} blows up in finite time.
\end{theorem}
\begin{exam}\label{Exam-01}
We here propose two typical classes of nonlinearities verifying the critical exponent $p_c(n)=1+2/n$ according to \eqref{GESDS-Condition} and \eqref{Blow-up-Condition-01}. Precisely, the global in time existence result in Theorem \ref{Thm-GESDS} holds when $p>p_c(n)$, whereas the blow-up result in Theorem \ref{Thm-Blow-up-01} holds when $1<p\leqslant p_c(n)$.
\begin{description}
	\item[Class I.1.] If one of the nonlinearities is vanishing, for example,
	\begin{align*}
	\begin{cases}
		f_1^N=|u_t^N|^{p}+|\theta^N|^{p}\\
		f_2^N=0
	\end{cases} \mbox{or}\ \ 	\begin{cases}
	f_1^N=0\\
	f_2^N=|\Delta u^N|^{p}+|u_t^N|^{p}+|\theta^N|^{p}
\end{cases}
\end{align*}
 then the critical exponent is given by $p=p_{\mathrm{c}}(n)$.
 \item[Class I.2.] If the power exponents are equal to $p$, for example,
 \begin{align*}
 \begin{cases}
  f_1^N=|\Delta u^N|^{p}+|u_t^N|^{p}\\
 f_2^N=|\Delta u^N|^{p}+|\theta^N|^{p}
 \end{cases}
\mbox{or}\ \
 \begin{cases}
	f_1^N=|u_t^N|^{p}+|\theta^N|^{p}\\
	f_2^N=|\Delta u^N|^{p}
\end{cases}
 \end{align*}
 then the critical exponent is given by $p=p_{\mathrm{c}}(n)$.
\end{description}
By taking $f_1^N=|u_t^N|^p$ and $f_2^N=0$, its critical exponent exactly coincides with the one for \eqref{Semilinear-Sigma} with the derivative-type nonlinearity $|u_t^N|^p$, i.e. the second situation in \eqref{Critical-sigma}.
\end{exam}

For another, concerning $c_{1,1}+c_{1,2}>0$, there exists at least a coefficient of  $|u^N|^{p_j}$ being non-vanishing. We thus may introduce another notation of weak solution to the semilinear thermoelastic plate systems \eqref{General-Semilinear-TEP}.
\begin{defn}\label{Defn-weak-02}
	 	A pair of functions $(u^N,\theta^N)$ is said to be a global in time weak solution to the semilinear Cauchy problem \eqref{General-Semilinear-TEP} with the nonlinearities \eqref{nonlinearity} carrying $c_{1,1}+c_{1,2}>0$ if
	 \begin{align*}
	 	\begin{cases}
	 		u^N\in L^{p_j}_{\mathrm{loc}}([0,T)\times \mb{R}^n)&\mbox{if}\ \ c_{1,j}>0,\\
	 		|\nabla |^{\sigma}u^N\in L^{p_j}_{\mathrm{loc}}([0,T)\times \mb{R}^n)&\mbox{if}\ \ c_{2,j}>0,\\
	 		u^N_t\in L^{p_j}_{\mathrm{loc}}([0,T)\times \mb{R}^n)&\mbox{if}\ \ c_{3,j}>0,\\
	 		\theta^N\in L^{p_j}_{\mathrm{loc}}([0,T)\times \mb{R}^n)&\mbox{if}\ \ c_{4,j}>0,		
	 	\end{cases}
	 \end{align*}
	 and	 the following integral equalities:
	\begin{align*}
		&\int_0^{+\infty}\int_{\mb{R}^n}f_1^N(t,x)\Psi_1(t,x)\,\mathrm{d}x\,\mathrm{d}t+\int_{\mb{R}^n}\big(-u^N_0(x)\,\partial_t\Psi_1(0,x)+u_1^N(x)\Psi_1(0,x)\big)\,\mathrm{d}x\\
		&=\int_0^{+\infty}\int_{\mb{R}^n}u^N(t,x)\big(\partial_{t}^2\Psi_1(t,x)+\Delta^2\Psi_1(t,x)\big)\,\mathrm{d}x\,\mathrm{d}t+\int_0^{+\infty}\int_{\mb{R}^n}\theta^N(t,x)\Delta \Psi_1(t,x)\,\mathrm{d}x\,\mathrm{d}t
	\end{align*}
as well as
	\begin{align*}
	&\int_0^{+\infty}\int_{\mb{R}^n}f_2^N(t,x)\Psi_2(t,x)\,\mathrm{d}x\,\mathrm{d}t+\int_{\mb{R}^n}\big(-u_0^N(x)\Delta\Psi_2(0,x)+\theta_0^N(x)\Psi_2(0,x)\big)\,\mathrm{d}x\\
	&=\int_0^{+\infty}\int_{\mb{R}^n}\theta^N(t,x)\big(-\partial_t\Psi_2(t,x)-\Delta \Psi_2(t,x)\big)\,\mathrm{d}x\,\mathrm{d}t+\int_0^{+\infty}\int_{\mb{R}^n}u^N(t,x)\Delta\partial_t \Psi_2(t,x)\,\mathrm{d}x\,\mathrm{d}t
\end{align*}
	hold for any test functions $\Psi_1,\Psi_2\in\ml{C}_0^{\infty}([0,+\infty)\times\mb{R}^n)$.
\end{defn}
\begin{theorem}\label{Thm-Blow-up-02}
	Let us consider the nonlinearities \eqref{nonlinearity} with $\sigma\in\{1,2\}$ satisfying
	\begin{itemize}
		\item $c_{d,j}\geqslant0$ for any $d\in\{1,\dots,4\}$ and $j\in\{1,2\}$;
		\item $c_{1,1}+c_{1,2}>0$,  $c_{4,1}+c_{4,2}>0$;
	\end{itemize}
	and their powers 
\begin{align}\label{Blow-up-Condition-02}
	1<p_j\leqslant \begin{cases}
		\displaystyle{1+\frac{4}{(n-2)_+}}&\mbox{if}\ \ c_{1,j}>0,\\[1em]
		\displaystyle{1+\frac{2}{n}}&\mbox{if}\ \ c_{4,j}>0,
	\end{cases}
\end{align}
	for $j\in\{1,2\}$. Suppose that the initial data $u_0^N=0$ and $u_1^N+\theta_0^N\in L^1$ verifying
	\begin{align*}
		\int_{\mb{R}^n}\big(u_1^N(x)+\theta_0^N(x)\big)\,\mathrm{d}x>0.
	\end{align*}
	Then, every local in time weak solution in the sense of Definition \ref{Defn-weak-02} to the general semilinear thermoelastic plate systems \eqref{General-Semilinear-TEP} blows up in finite time.
\end{theorem}

\begin{exam}\label{Exam-02}
	We here propose two typical classes of nonlinearities verifying the critical exponent $\bar{p}_c(n)=1+4/(n-2)_+$ according to \eqref{GESDS-Condition} and \eqref{Blow-up-Condition-02}. Precisely, the global in time existence result in Theorem \ref{Thm-GESDS} holds when $p>\bar{p}_c(n)$, whereas the blow-up result in Theorem \ref{Thm-Blow-up-02} holds when $1<p\leqslant \bar{p}_c(n)$. Recall here that $\sigma\in\{1,2\}$.
	\begin{description}
		\item[Class II.1.] If one of the nonlinearities is vanishing, for example,
		\begin{align*}
			\begin{cases}
				f_1^N=|u^N|^{p}+|\theta^N|^{p}\\
				f_2^N=0
			\end{cases} \mbox{or}\ \ 	\begin{cases}
				f_1^N=0\\
				f_2^N=|u^N|^{p}+|\,|\nabla|^{\sigma} u^N|^{p}+|u_t^N|^{p}+|\theta^N|^{p}
			\end{cases}
		\end{align*}
		then the critical exponent is given by $p=\bar{p}_{\mathrm{c}}(n)$.
		\item[Class II.2.] If the power exponents are suitably proportional related to $p$, for example,
		\begin{align*}
			\begin{cases}
				f_1^N=|u|^p+|\,|\nabla|^{\sigma} u^N|^{p}+|u_t^N|^{p}\\
				f_2^N=|\,|\nabla|^{\sigma} u^N|^{\frac{(n-2)_+}{n}p}+|\theta^N|^{\frac{(n-2)_+}{n}p}
			\end{cases}
			\mbox{or}\ \
			\begin{cases}
				f_1^N=|u^N|^{p}+|\theta^N|^{p}\\
				f_2^N=|u^N|^p+|\,|\nabla|^{\sigma} u^N|^{p}+|u^N_t|^p
			\end{cases}
		\end{align*}
		then the critical exponent is given by $p=\bar{p}_{\mathrm{c}}(n)$.
	\end{description}
By taking $f_1^N=|u^N|^p$ and $f_2^N=0$, its critical exponent exactly coincides with the one for \eqref{Semilinear-Sigma} with the power-type nonlinearity $|u^N|^p$, i.e. the first situation in \eqref{Critical-sigma}.
\end{exam}

\section{$(L^q\cap L^m)-L^q$ estimates for the linear thermoelastic plate system}\setcounter{equation}{0}\label{Section-Linear}
\hspace{5mm}This section contributes to derive some new $(L^q\cap L^m)-L^q$ estimates with any $1\leqslant m\leqslant q\leqslant +\infty$ of the vertical displacement $u$ and the temperature difference $\theta$ themselves and their time/spatial derivatives, respectively, to the linear thermoelastic plate system \eqref{Linear-TEP}. Our study of qualitative properties of the solutions $u,\theta$ not only is significant to understand underlying physical phenomena but also devotes to study global in time behavior for its corresponding nonlinear model \eqref{General-Semilinear-TEP}. Although the $L^q$ theory for the linear thermoelastic plate system \eqref{Linear-TEP} has been deeply developed by \cite{Denk-Racke=2006} in the whole space when $q\in[2,+\infty]$, and by \cite{Denk-Racke-Shibata=2009,Naito-Shibata=2009} in some bounded or exterior domains when $q\in(1,+\infty)$, qualitative properties of solutions (particularly the solution $u$ itself) in the general $L^q$ framework with $q\in[1,+\infty]$ are still not fixed.

We stress that our derivation of estimates of $u(t,\cdot)$ itself in the general $L^q$ framework is not simply an application of the pointwise estimates of the energy term $\mb{E}(t,\cdot)$ in \cite{Denk-Racke=2006,Said-Houari=2013,Racke-Ueda=2016}. \begin{itemize}
	\item For one thing, with the aid of sharp quantitative estimate in \cite[Inequality (2.3)]{Racke-Ueda=2016}, that is
	\begin{align}\label{Racke-Pointwise}
		|\xi|^2|\widehat{u}|\leqslant\frac{13}{3}\,\mathrm{e}^{-\frac{1}{52}|\xi|^2t}\left(|\xi|^2|\widehat{u}_0|+|\widehat{u}_1|+|\widehat{\theta}_0|\right),
	\end{align}
	one just can deduce the decay estimate for the energy term $\Delta u(t,\cdot)$ in the $L^q$ norm with $q\in[2,+\infty]$ by the Hausdorff-Young inequality as follows:
	\begin{align}\label{Estimate-Delta-u}
		\|\Delta u(t,\cdot)\|_{L^q}\lesssim \|\,|\xi|^2\widehat{u}(t,\xi)\|_{L^{q'}}&\lesssim \left\|\,\mathrm{e}^{-\frac{1}{52}|\xi|^2t}\left(|\xi|^2|\widehat{u}_0(\xi)|+|\widehat{u}_1(\xi)|+|\widehat{\theta}_0(\xi)|\right)\right\|_{L^{q'}}\ \ \mbox{with}\ \ q'\in[1,2]\notag\\
		&\lesssim \left\|\,\mathrm{e}^{-c|\xi|^2t}\right\|_{L^{\frac{m'q'}{m'-q'}}}\left\||\xi|^2|\widehat{u}_0(\xi)|+|\widehat{u}_1(\xi)|+|\widehat{\theta}_0(\xi)|\right\|_{L^{m'}}\notag\\
		&\lesssim t^{-\frac{n}{2}\left(\frac{1}{m}-\frac{1}{q}\right)}\left(\|u_0\|_{\dot{H}^2_m}+\|u_1\|_{L^m}+\|\theta_0\|_{L^m}\right)
	\end{align}
	for $1\leqslant m\leqslant 2\leqslant q\leqslant +\infty$.  In the above deduction with $m<q$,  the next estimate was used: 
	\begin{align*}
		\left\|\,\mathrm{e}^{-c|\xi|^2t}\right\|_{L^{\frac{m'q'}{m'-q'}}}=\left\|\,\mathrm{e}^{-c|\xi|^2t}\right\|_{L^{\frac{mq}{q-m}}}&\lesssim\left(\,\int_0^{+\infty}\mathrm{e}^{-\frac{mqc}{q-m}|\xi|^2t}\,\mathrm{d}|\xi|^n\right)^{\frac{q-m}{mq}}\lesssim t^{-\frac{n}{2}\left(\frac{1}{m}-\frac{1}{q}\right)}.
	\end{align*}
This is a simple generalization of \cite[Theorem 4.2]{Denk-Racke=2006}, where we do not need the conjugate relation between $m$ and $q$ in \eqref{Estimate-Delta-u}.
	However, \eqref{Estimate-Delta-u} and \cite[Theorem 4.2]{Denk-Racke=2006} cannot show any estimate for the case $q\in[1,2)$ or $m\in(2,+\infty]$.
	\item For another, by using the Hardy-Littlewood-Sobolev inequality and the estimate \eqref{Estimate-Delta-u} for the energy term $\Delta u(t,\cdot)$, one may find
	\begin{align*}
		\|u(t,\cdot)\|_{L^q}=\|\,|D|^{-2}\Delta u(t,\cdot)\|_{L^q}&\lesssim\|\Delta u(t,\cdot)\|_{L^{\frac{nq}{2q+n}}}\ \ \mbox{with}\ \ 1<\frac{nq}{2q+n}<\frac{n}{2}\ \ \mbox{and}\ \ q<+\infty\\
		&\lesssim t^{1-\frac{n}{2}\left(\frac{1}{m}-\frac{1}{q}\right)}\left(\|u_0\|_{\dot{H}^2_m}+\|u_1\|_{L^m}+\|\theta_0\|_{L^m}\right)
	\end{align*}
	with the technical restrictions
	\begin{align*}
		1\leqslant m\leqslant 2\leqslant q\leqslant +\infty,\ \ \frac{n}{n-2}<q<+\infty,\ \ n\geqslant 3.
	\end{align*}
	In other words, the questions on lower dimensions and other value of $q,m$ are still open. Additionally, the data space $\dot{H}^2_m$ for $u_0$ seems unusual due to the desired estimate for $u(t,\cdot)$ itself in the $L^q$ norm only.
\end{itemize} 

 To overcome these difficulties, motivated by the recent paper \cite{Chen-Ikehata=2023} we recover some oscillations to compensate strong singularities via the explicit representation of solutions in the Fourier space rather than an application of the Hardy-Littlewood-Sobolev inequality. Then, according to $L^r$ estimates for oscillating integrals instead of the Hausdorff-Young inequality, we derive the sharp estimates for $|D|^s\partial_t^ku(t,\cdot)$ and $|D|^s\theta(t,\cdot)$ with any $s\in[0,+\infty)$ and $k\in\{0,1\}$ in the general $L^q$ norm for any $q\in[1,+\infty]$. As byproducts, the $L^q$ well-posedness and the generalized Gevrey (ultra-analytic) smoothing effect of solutions are deduced.

\subsection{Explicit representation of solutions in the Fourier space}
\hspace{5mm}Let us follow the reduction procedure in \cite[Section 3]{Chen-Ikehata=2023}. Precisely, by acting the heat operator $\partial_t-\Delta$ on \eqref{Linear-TEP}$_1$ and summing up the resultant with \eqref{Linear-TEP}$_2$, the vertical displacement $u$ satisfies the higher order evolution equation
\begin{align}\label{Eq-third-order}
	\begin{cases}
		u_{ttt}-\Delta u_{tt}+2\Delta^2 u_t-\Delta^3 u=0,&x\in\mb{R}^n,\ t>0,\\
		(u,u_t,u_{tt})(0,x)=(u_0,u_1,-\Delta^2 u_0-\Delta\theta_0)(x),&x\in\mb{R}^n,
	\end{cases}
\end{align}
moreover, the temperature difference $\theta$ can be represented according to
\begin{align}\label{Relation-theta-01}
\theta=(-\Delta)^{-1}u_{tt}-\Delta u.
\end{align}
The characteristic equation to \eqref{Eq-third-order} is addressed by $\lambda^3+|\xi|^2\lambda^2+2|\xi|^4\lambda+|\xi|^6=0$ whose roots are
\begin{align*}
	\lambda_1&=-a_0|\xi|^2:=-\frac{1+\alpha_-}{3}|\xi|^2,\\
	\lambda_{2,3}&=-a_1|\xi|^2\mp ia_2|\xi|^2:=-\frac{2-\alpha_-}{6}|\xi|^2\mp i\frac{\sqrt{3}\,\alpha_+}{6}|\xi|^2,
\end{align*}
where two numbers $\alpha_{\pm}$ are defined by
\begin{align*}
	\alpha_{\pm}:=\sqrt[3]{\frac{1}{2}(3\sqrt{69}+11)}\pm\sqrt[3]{\frac{1}{2}(3\sqrt{69}-11)}.
\end{align*}
These parameters also coincide with those in \cite[Subsection 2.2]{Racke-Ueda=2016}. Note that $a_0\approx 0.57$, $a_1\approx 0.22$ and $a_2\approx 1.31$. For the sake of simplicity, we may denote $a_3:=(a_0-a_1)^2+a_2^2>0$.

Thanks to these pairwise distinct characteristic roots, the solutions in the Fourier space, i.e. $\widehat{u}=\widehat{u}(t,\xi)$ and $\widehat{\theta}=\widehat{\theta}(t,\xi)$, are given by
\begin{align}\label{Rep-u-01}
	\widehat{u}&=\widehat{K}_0(t,|\xi|)\,\widehat{u}_0+\widehat{K}_1(t,|\xi|)\,\widehat{u}_1+\widehat{K}_2(t,|\xi|)\left(-|\xi|^4\widehat{u}_0+|\xi|^2\widehat{\theta}_0\right)\notag\\
	&=\big(\widehat{K}_0(t,|\xi|)-|\xi|^4\widehat{K}_2(t,|\xi|)\big)\,\widehat{u}_0+\widehat{K}_1(t,|\xi|)\,\widehat{u}_1+|\xi|^2\widehat{K}_2(t,|\xi|)\,\widehat{\theta}_0\notag\\
	&=:\widehat{K}_0^u(t,|\xi|)\,\widehat{u}_0+\widehat{K}_1^u(t,|\xi|)\,\widehat{u}_1+\widehat{K}_2^u(t,|\xi|)\,\widehat{\theta}_0
\end{align}
and, from \eqref{Relation-theta-01},
\begin{align}\label{Rep-theta-01}
	\widehat{\theta}&=\left((|\xi|^{-2}\partial_t^2+|\xi|^2)\widehat{K}_0(t,|\xi|)-(|\xi|^2\partial_t^2+|\xi|^6)\widehat{K}_2(t,|\xi|)\right)\widehat{u}_0\notag\\
	&\quad+ (|\xi|^{-2}\partial_t^2+|\xi|^2)\widehat{K}_1(t,|\xi|)\,\widehat{u}_1+(\partial_t^2+|\xi|^4)\widehat{K}_2(t,|\xi|)\,\widehat{\theta}_0\notag\\
	&=:\widehat{K}_0^{\theta}(t,|\xi|)\,\widehat{u}_{0}+\widehat{K}_1^{\theta}(t,|\xi|)\,\widehat{u}_1+\widehat{K}_2^{\theta}(t,|\xi|)\,\widehat{\theta}_0.
\end{align}
In the above, we re-combine the multipliers in an appropriate way (due to some compensations that will be shown later)
\begin{align*}
	\widehat{K}_0(t,|\xi|)&=\frac{a_1^2+a_2^2}{a_3}\,\mathrm{e}^{-a_0|\xi|^2t}+\frac{a_0^2-2a_0a_1}{2a_3}\,\widehat{G}_+(t,|\xi|)+\frac{a_0(a_0a_1-a_1^2+a_2^2)}{2ia_2a_3}\,\widehat{G}_-(t,|\xi|),\\
	\widehat{K}_1(t,|\xi|)&=\frac{a_1}{a_3|\xi|^2}\left( 2\,\mathrm{e}^{-a_0|\xi|^2t}-\widehat{G}_+(t,|\xi|)\right)+\frac{a_0^2+a_2^2-a_1^2}{2ia_2a_3|\xi|^2}\,\widehat{G}_-(t,|\xi|),\\
	\widehat{K}_2(t,|\xi|)&=\frac{1}{2a_3|\xi|^4}\left(2\,\mathrm{e}^{-a_0|\xi|^2t}-\widehat{G}_{+}(t,|\xi|)\right)+\frac{a_0-a_1}{2ia_2a_3|\xi|^4}\,\widehat{G}_-(t,|\xi|),
\end{align*}
with the diffusion plate factors (see \cite{Ikehata=2021} with $\sigma=2$)
\begin{align*}
	\widehat{G}_{\pm}(t,|\xi|):=\mathrm{e}^{(-a_1+ia_2)|\xi|^2t}\pm\mathrm{e}^{(-a_1-ia_2)|\xi|^2t}=\begin{cases}
	2\cos(a_2|\xi|^2t)\,\mathrm{e}^{-a_1|\xi|^2t}&\mbox{if} \ \ \mbox{``$+$''},\\
	2i\sin(a_2|\xi|^2t)\,\mathrm{e}^{-a_1|\xi|^2t}&\mbox{if} \ \ \mbox{``$-$''}.
	\end{cases}
\end{align*}
Although the representations of $\widehat{u},\widehat{\theta}$ have been obtained by \cite[Formulas (2.12) and (2.14)]{Naito-Shibata=2009}, our re-arrangement may contribute to the endpoint estimates in the $L^q$ norm with $q\in\{1,+\infty\}$.
\begin{remark}\label{Rem-Pointwise}
From the solutions' formulas \eqref{Rep-u-01} and \eqref{Rep-theta-01} in the Fourier space, one easily gets
\begin{align*}
|\widehat{u}|&\lesssim \mathrm{e}^{-c|\xi|^2t}\left(|\widehat{u}_0|+|\xi|^{-2}|\widehat{u}_1|+|\xi|^{-2}|\widehat{\theta}_0|\right),\\
|\widehat{\theta}|&\lesssim \mathrm{e}^{-c|\xi|^2t}\left(|\xi|^2|\widehat{u}_0|+|\widehat{u}_1|+|\widehat{\theta}_0|\right),
\end{align*}
which imply regularity-gain-type dissipative mechanism (proposed by \cite[Formula (1.4)]{Kawashima-Shibata-Xu=2022} for the linearized compressible Navier-Stokes-Korteweg system) of the linear thermoelastic system \eqref{Linear-TEP} because of the heat-type kernel $\mathrm{e}^{-c|\xi|^2t}$ for all frequencies $\xi\in\mb{R}^n$.
\end{remark}

\subsection{Sharp $L^r$ estimates for the crucial Fourier multipliers}
\hspace{5mm}One notices that the kernels $\widehat{K}_j(t,|\xi|)$ consist of the diffusion plate factors $\widehat{G}_{\pm}(t,|\xi|)$ with the singularities $|\xi|^{-2}$ or $|\xi|^{-4}$, and the Gaussian kernel $\mathrm{e}^{-a_0|\xi|^2t}$, in the sense of linear combinations. As preparations of estimating the solutions in the general $L^q$ norm, we next propose some sharp $L^r$ estimates for these crucial Fourier multipliers in the $L^r$ norm with $r\in[1,+\infty]$.

Let us first recall \cite[Propositions 5, 12]{Pham-Kainae-Reissig=2015} and \cite[Subsection 24.2]{Ebert-Reissig=2018}, which stated the next lemma for the Gaussian kernel even with oscillations. It contributes to estimating the kernel $K_0(t,|D|)$.
\begin{lemma}\label{Lemma-Kernel}
	Let $c_0>0$ and $c_1\neq 0$. Then, concerning $s\in[0,+\infty)$ and $r\in[1,+\infty]$, the following decay estimates:
	\begin{align*}
		\left\|\ml{F}^{-1}_{\xi\to x}\left(|\xi|^{s}\,\mathrm{e}^{-c_0|\xi|^2t}g_0(c_1|\xi|^2t)\right)\right\|_{L^r}\lesssim t^{-\frac{n}{2}\left(1-\frac{1}{r}\right)-\frac{s}{2}}
	\end{align*}
	hold for any $t>0$, with $g_0(\tau)\in\{1,\sin\tau,\cos\tau\}$.
\end{lemma}

For the other two kernels, there are two cancellations: between the Gaussian kernel and $\widehat{G}_+(t,|\xi|)$; in $\widehat{G}_-(t,|\xi|)$ itself, which produce the additional factor $|\xi|^2$ to overcome the strong singularity $|\xi|^{-2}$ in the kernels $K_1(t,|D|)$ and $|D|^2K_2(t,|D|)$. Thus, the estimate \eqref{Racke-Pointwise} can be improved in the sense of recovering some oscillations.
\begin{lemma}\label{Lemma-Kernel-2}
	Let $s\in[0,+\infty)$ and $k\in\mb{N}_0$. Then, concerning $r\in[1,+\infty]$, the following sharp estimates:
	\begin{align}
		\left\||D|^s\partial_t^k\,\ml{F}^{-1}_{\xi\to x}\left( \frac{1}{|\xi|^2}\left( 2\,\mathrm{e}^{-a_0|\xi|^2t}-\widehat{G}_+(t,|\xi|)\right)\right)\right\|_{L^r}&\lesssim t^{-\frac{n}{2}\left(1-\frac{1}{r}\right)-\frac{s}{2}-k+1},\label{Est-02}\\
		\left\||D|^s\partial_t^k\,\ml{F}^{-1}_{\xi\to x}\left( \frac{1}{2i|\xi|^2}\,\widehat{G}_-(t,|\xi|)\right)\right\|_{L^r}&\lesssim t^{-\frac{n}{2}\left(1-\frac{1}{r}\right)-\frac{s}{2}-k+1},\label{Est-03}
	\end{align}
	hold for any $t>0$.
\end{lemma}

\begin{proof}
To begin with the proof, let us rewrite the multiplier in an appropriate way
\begin{align*}
	\frac{1}{|\xi|^2}\left( 2\,\mathrm{e}^{-a_0|\xi|^2t}-\widehat{G}_+(t,|\xi|)\right)&=\sum\limits_{\pm}\frac{1}{|\xi|^2}\left(\mathrm{e}^{-a_0|\xi|^2t}-\mathrm{e}^{(-a_1\pm ia_2)|\xi|^2t}\right)\\
	&=-t\,\mathrm{e}^{-a_0|\xi|^2t}\sum\limits_{\pm}(a_0-a_1\pm ia_2)\int_0^1\mathrm{e}^{(a_0-a_1\pm ia_2)|\xi|^2t\tau}\,\mathrm{d}\tau
\end{align*}
and, for $k\in\mb{N}_+$,
\begin{align*}
	\frac{1}{|\xi|^2}\,\partial_t^k\left( 2\,\mathrm{e}^{-a_0|\xi|^2t}-\widehat{G}_+(t,|\xi|)\right)=|\xi|^{2k-2}\sum\limits_{\pm}\left((-a_0)^k\,\mathrm{e}^{-a_0|\xi|^2t}-(-a_1\pm ia_2)^k\,\mathrm{e}^{(-a_1\pm ia_2)|\xi|^2t}\right).
\end{align*}
Indeed, Lemma \ref{Lemma-Kernel} implies
\begin{align}\label{Est-04}
		\left\|\ml{F}^{-1}_{\xi\to x}\left(|\xi|^{s}\,\mathrm{e}^{(-c_0+ic_1)|\xi|^2t}\right)\right\|_{L^r}\lesssim t^{-\frac{n}{2}\left(1-\frac{1}{r}\right)-\frac{s}{2}}
\end{align}
for any $t>0$, where $c_0>0$, $c_1\in\mb{R}$, $s\in[0,+\infty)$ and $r\in[1,+\infty]$. As a consequence,
\begin{align*}
	&\left\||D|^s\,\ml{F}^{-1}_{\xi\to x}\left( \frac{1}{|\xi|^2}\left( 2\,\mathrm{e}^{-a_0|\xi|^2t}-\widehat{G}_+(t,|\xi|)\right)\right)\right\|_{L^r}\\
	&\lesssim t\sum\limits_{\pm}\left\|\ml{F}_{\xi\to x}^{-1}\left(|\xi|^{s}\, \mathrm{e}^{-a_0|\xi|^2t}\int_0^1\mathrm{e}^{(a_0-a_1\pm ia_2)|\xi|^2t\tau}\,\mathrm{d}\tau \right)\right\|_{L^r}\\
	&\lesssim t\sum\limits_{\pm}\int_0^1\left\|\ml{F}_{\xi\to x}^{-1}\left(|\xi|^s\,\mathrm{e}^{[-a_0+(a_0-a_1\pm ia_2)\tau]|\xi|^2t}\right)\right\|_{L^r}\mathrm{d}\tau\\
	&\lesssim t^{-\frac{n}{2}\left(1-\frac{1}{r}\right)-\frac{s}{2}+1}
\end{align*}
holds for $s\in[0,+\infty)$ and $r\in[1,+\infty]$, due to the fact that $a_0-(a_0-a_1)\tau=(1-\tau)a_0+a_1\tau>0$ for any $\tau\in[0,1]$. Similarly,
\begin{align*}
	&\left\||D|^s\partial_t^k\,\ml{F}^{-1}_{\xi\to x}\left( \frac{1}{|\xi|^2}\left( 2\,\mathrm{e}^{-a_0|\xi|^2t}-\widehat{G}_+(t,|\xi|)\right)\right)\right\|_{L^r}\\
	&\lesssim\left\|\ml{F}^{-1}_{\xi\to x}\left(|\xi|^{s+2k-2}\,\mathrm{e}^{-a_0|\xi|^2t}\right)\right\|_{L^r}+\sum\limits_{\pm}\left\|\ml{F}^{-1}_{\xi\to x}\left(|\xi|^{s+2k-2}\,\mathrm{e}^{(-a_1\pm ia_2)|\xi|^2t}\right)\right\|_{L^r}\\
	&\lesssim t^{-\frac{n}{2}\left(1-\frac{1}{r}\right)-\frac{s}{2}-k+1}
\end{align*}
holds for $s+2k-2\geqslant0$ if $k\in\mb{N}_+$. The last two estimates conclude \eqref{Est-02} immediately.

For another, thanks to
\begin{align*}
	\frac{1}{2i|\xi|^2}\,\widehat{G}_-(t,|\xi|)&=\frac{1}{2i|\xi|^2}\,\mathrm{e}^{(-a_1+ia_2)|\xi|^2t}\left(1-\mathrm{e}^{-2ia_2|\xi|^2t}\right)\\
	&=a_2\,t\,\mathrm{e}^{(-a_1+ia_2)|\xi|^2t}\int_0^1\mathrm{e}^{-2ia_2|\xi|^2t\tau}\,\mathrm{d}\tau
\end{align*}
as well as, for $k\in\mb{N}_+$,
\begin{align*}
	\frac{1}{2i|\xi|^2}\,\partial_t^k\widehat{G}_-(t,|\xi|)=\frac{|\xi|^{2k-2}}{2i}\left((-a_1+ia_2)^k\,\mathrm{e}^{(-a_1+ia_2)|\xi|^2t}-(-a_1-ia_2)^k\,\mathrm{e}^{(-a_1-ia_2)|\xi|^2t}\right),
\end{align*}
by using \eqref{Est-04} we are able to derive
\begin{align*}
	\left\||D|^s\partial_t^k\,\ml{F}^{-1}_{\xi\to x}\left( \frac{1}{2i|\xi|^2}\,\widehat{G}_-(t,|\xi|)\right)\right\|_{L^r}
	&\lesssim\begin{cases}
	\displaystyle{t\int_0^1\left\|\ml{F}_{\xi\to x}^{-1}\left(|\xi|^s\,\mathrm{e}^{[-a_1+ia_2(1-2\tau)]|\xi|^2t}\right)\right\|_{L^r}\mathrm{d}\tau}&\mbox{if}\ \ k=0\\[1em]
		\displaystyle{\sum\limits_{\pm}\left\|\ml{F}_{\xi\to x}^{-1}\left(|\xi|^{s+2k-2}\,\mathrm{e}^{(-a_1\pm ia_2)|\xi|^2t}\right)\right\|_{L^r}}&\mbox{if}\ \ k\in\mb{N}_+
	\end{cases}\\
&\lesssim t^{-\frac{n}{2}\left(1-\frac{1}{r}\right)-\frac{s}{2}-k+1}
\end{align*}
for $s\geqslant0$ if $k=0$ and $s+2k-2\geqslant 0$ if $k\in\mb{N}_+$, which demonstrates \eqref{Est-03}.
\end{proof}
\begin{remark}
In \cite[Proposition 3]{Chen-Ikehata=2023}, the following large time optimal estimates:
	\begin{align*}
	\left\|\ml{F}^{-1}_{\xi\to x}\left( \frac{1}{|\xi|^2}\left( 2\,\mathrm{e}^{-a_0|\xi|^2t}-\widehat{G}_+(t,|\xi|)\right)\right)\right\|_{L^2},\ \ 
	\left\|\ml{F}^{-1}_{\xi\to x}\left( \frac{1}{2i|\xi|^2}\,\widehat{G}_-(t,|\xi|)\right)\right\|_{L^2}\approx t^{1-\frac{n}{4}}
\end{align*}
have been proved by the elementary analysis in the Fourier space thanks to the validity of Plancherel theorem in $L^2$, which exactly coincide with our result in Lemma \ref{Lemma-Kernel-2} with $s=k=0$ and $r=2$. For this reason, we expect our estimates are sharp.
\end{remark}

\subsection{$(L^q\cap L^m)-L^q$ estimates for the solutions}
\hspace{5mm}We now state the main result on the linear thermoelastic plate system \eqref{Linear-TEP}.
\begin{theorem}\label{Thm-Linear-u}
Let $s\in[0,+\infty)$ and $k\in\{0,1\}$. Then, the vertical displacement $u$ and the temperature difference $\theta$, respectively, fulfill the following $(L^q\cap L^m)-L^q$ estimates:
\begin{align}\label{Est-01}
\|\,|D|^s\partial_t^k u(t,\cdot)\|_{L^q}&\lesssim (1+t)^{-\frac{n}{2}\left(\frac{1}{m}-\frac{1}{q}\right)-\frac{s}{2}-k}\|u_0\|_{\dot{H}^{s+2k}_q\cap L^m}+(1+t)^{-\frac{n}{2}\left(\frac{1}{m}-\frac{1}{q}\right)-\frac{s}{2}-k+1}\|u_1\|_{\dot{H}^{\max\{s+2k-2,0\}}_q\cap L^m}\notag\\
&\quad\ +(1+t)^{-\frac{n}{2}\left(\frac{1}{m}-\frac{1}{q}\right)-\frac{s}{2}-k+1}\|\theta_0\|_{\dot{H}^{\max\{s+2k-2,0\}}_q\cap L^m}
\end{align}
and 	
\begin{align*}
	\|\,|D|^s\theta(t,\cdot)\|_{L^q}&\lesssim (1+t)^{-\frac{n}{2}\left(\frac{1}{m}-\frac{1}{q}\right)-\frac{s}{2}-1}\|u_0\|_{\dot{H}^{s+2}_q\cap L^m}+(1+t)^{-\frac{n}{2}\left(\frac{1}{m}-\frac{1}{q}\right)-\frac{s}{2}}\|u_1\|_{\dot{H}^{s}_q\cap L^m}\\
	&\quad+(1+t)^{-\frac{n}{2}\left(\frac{1}{m}-\frac{1}{q}\right)-\frac{s}{2}}\|\theta_0\|_{\dot{H}^{s}_q\cap L^m}
\end{align*}
for any $t\geqslant0$ and any $1\leqslant m\leqslant q\leqslant +\infty$. Particularly, when $m=q$, they are the so-called $L^q-L^q$ estimates.
\end{theorem}

\begin{proof}
The explicit representation of $\widehat{K}_0^u(t,|\xi|)$ associated with \eqref{Est-04} shows
\begin{align*}
	&\left\| |D|^s\partial_t^k\,\ml{F}^{-1}_{\xi\to x}\left(\widehat{K}_0^u(t,|\xi|)\right)u_0(\cdot)\right\|_{L^q}\\
	&\lesssim \left(\,\left\| |D|^s\partial_t^k\,\ml{F}^{-1}_{\xi\to x}\left(\widehat{K}_0(t,|\xi|)\right)\right\|_{L^r}+\left\| |D|^{s}\partial_t^k\,\ml{F}^{-1}_{\xi\to x}\left(|\xi|^4\widehat{K}_2(t,|\xi|)\right)\right\|_{L^r}\right)\|u_0\|_{L^m}\\
	&\lesssim \left(\left\||D|^s\partial_t^k\,\ml{F}_{\xi\to x}^{-1}\left(\mathrm{e}^{-a_0|\xi|^2t}\right)\right\|_{L^r}+\sum\limits_{\pm}\left\||D|^s\partial_t^k\,\ml{F}_{\xi\to x}^{-1}\left(\widehat{G}_{\pm}(t,|\xi|)\right)\right\|_{L^r}\right)\|u_0\|_{L^m}\\
	&\lesssim \left\|\ml{F}^{-1}_{\xi\to x}\left(|\xi|^{s+2k}\,\mathrm{e}^{(-c_0+ic_1)|\xi|^2t}\right)\right\|_{L^r}\|u_0\|_{L^m}\ \ \mbox{with suitable}\ \ c_0>0,\ c_1\in\mb{R}\\
	&\lesssim t^{-\frac{n}{2}\left(1-\frac{1}{r}\right)-\frac{s}{2}-k}\|u_0\|_{L^m}
\end{align*}
for $t>1$, $s\in[0,+\infty)$ and $k\in\{0,1\}$, where $1+1/q=1/r+1/m$ with $1\leqslant m\leqslant q\leqslant +\infty$ due to an application of Young's convolution inequality. Moreover, to avoid the singularity as $t\to0^+$ one may directly use the next bounded estimate:
\begin{align}
	&\left\| |D|^s\partial_t^k\,\ml{F}^{-1}_{\xi\to x}\left(\widehat{K}_0^u(t,|\xi|)\right)u_0(\cdot)\right\|_{L^q}\notag\\
	&\lesssim \left(\left\||D|^{-2k}\partial_t^k\,\ml{F}_{\xi\to x}^{-1}\left(\mathrm{e}^{-a_0|\xi|^2t}\right)\right\|_{L^1}+\sum\limits_{\pm}\left\||D|^{-2k}\partial_t^k\,\ml{F}_{\xi\to x}^{-1}\left(\widehat{G}_{\pm}(t,|\xi|)\right)\right\|_{L^1}\right)\|\,|D|^{s+2k}u_0\|_{L^q}\notag\\
	&\lesssim \left\|\ml{F}^{-1}_{\xi\to x}\left(\mathrm{e}^{(-c_0+ic_1)|\xi|^2t}\right)\right\|_{L^1}\|\,|D|^{s+2k}u_0\|_{L^q}\ \ \mbox{with suitable}\ \ c_0>0,\ c_1\in\mb{R}\notag\\
	&\lesssim \|u_0\|_{\dot{H}^{s+2k}_q}\label{Exist-01}
\end{align}
for $t\leqslant 1$, $s\in[0,+\infty)$ and $k\in\{0,1\}$. The summary of them claims
\begin{align*}
\left\| |D|^s\partial_t^k\,\ml{F}^{-1}_{\xi\to x}\left(\widehat{K}_0^u(t,|\xi|)\right)u_0(\cdot)\right\|_{L^q}\lesssim (1+t)^{-\frac{n}{2}\left(1-\frac{1}{r}\right)-\frac{s}{2}-k}\|u_0\|_{\dot{H}^{s+2k}_q\cap L^m}
\end{align*}
for any $t\geqslant 0$, $s\in[0,+\infty)$ and $k\in\{0,1\}$.

Next, by the same way as the above we use Lemma \ref{Lemma-Kernel-2} to derive
\begin{align*}
	&\left\| |D|^s\partial_t^k\,\ml{F}^{-1}_{\xi\to x}\left(\widehat{K}_1^u(t,|\xi|)\right)u_1(\cdot)\right\|_{L^q}\\
	&\lesssim \left(\,\left\||D|^s\partial_t^k\,\ml{F}^{-1}_{\xi\to x}\left( \frac{1}{|\xi|^2}\left( 2\,\mathrm{e}^{-a_0|\xi|^2t}-\widehat{G}_+(t,|\xi|)\right)\right)\right\|_{L^r}+\left\||D|^s\partial_t^k\,\ml{F}^{-1}_{\xi\to x}\left( \frac{1}{2i|\xi|^2}\,\widehat{G}_-(t,|\xi|)\right)\right\|_{L^r}\right)\|u_1\|_{L^m}\\
	&\lesssim t^{-\frac{n}{2}\left(1-\frac{1}{r}\right)-\frac{s}{2}-k+1}\|u_1\|_{L^m}
\end{align*}
for $t>1$, and
\begin{align}\label{Exist-02}
	&\left\| |D|^s\partial_t^k\,\ml{F}^{-1}_{\xi\to x}\left(\widehat{K}_1^u(t,|\xi|)\right)u_1(\cdot)\right\|_{L^q}\notag\\
	&\lesssim \left(\,\left\||D|^{\min\{2-2k,s\}}\partial_t^k\,\ml{F}^{-1}_{\xi\to x}\left( \frac{1}{|\xi|^2}\left( 2\,\mathrm{e}^{-a_0|\xi|^2t}-\widehat{G}_+(t,|\xi|)\right)\right)\right\|_{L^1}\right.\notag\\
	&\qquad\left.+\left\||D|^{\min\{2-2k,s\}}\partial_t^k\,\ml{F}^{-1}_{\xi\to x}\left( \frac{1}{2i|\xi|^2}\,\widehat{G}_-(t,|\xi|)\right)\right\|_{L^1}\right)\|\,|D|^{\max\{s+2k-2,0\}}u_1\|_{L^q}\notag\\
	&\lesssim t^{\frac{\max\{0,2-s-2k\}}{2}}\|u_1\|_{\dot{H}^{\max\{s+2k-2,0\}}_q} \lesssim\|u_1\|_{\dot{H}^{\max\{s+2k-2,0\}}_q} 
\end{align}
for $t\leqslant 1$, where $s\in[0,+\infty)$ and $k\in\{0,1\}$ guarantee $\min\{2-2k,s\}\geqslant0$. In the last estimate, we considered the equality
\begin{align*}
s=\max\{s+2k-2,0\}+\min\{2-2k,s\}.
\end{align*}
In other words,
\begin{align*}
\left\| |D|^s\partial_t^k\,\ml{F}^{-1}_{\xi\to x}\left(\widehat{K}_1^u(t,|\xi|)\right)u_1(\cdot)\right\|_{L^q}\lesssim (1+t)^{-\frac{n}{2}\left(1-\frac{1}{r}\right)-\frac{s}{2}-k+1}\|u_1\|_{\dot{H}^{\max\{s+2k-2,0\}}_q\cap L^m}
\end{align*}
holds for any $t\geqslant0$.

Thanks to $|\xi|^2\widehat{K}_2(t,|\xi|)\sim\widehat{K}_1(t,|\xi|)$, it is easy to check
\begin{align}\label{Exist-03}
	\left\| |D|^s\partial_t^k\,\ml{F}^{-1}_{\xi\to x}\left(\widehat{K}_2^u(t,|\xi|)\right)\theta_0(\cdot)\right\|_{L^q}\lesssim\begin{cases}
		t^{-\frac{n}{2}\left(1-\frac{1}{r}\right)-\frac{s}{2}-k+1}\|\theta_0\|_{L^m}&\mbox{for}\ \ t>1,\\
		\|\theta_0\|_{\dot{H}^{\max\{s+2k-2,0\}}_q} &\mbox{for}\ \ t\leqslant 1,
	\end{cases}
\end{align}
where $s\in[0,+\infty)$ and $k\in\{0,1\}$.

Finally, due to the fact that
\begin{align*}
 \partial_t^ku(t,\cdot)=\partial_t^k\left[\ml{F}^{-1}_{\xi\to x}\left(\widehat{K}_0^u(t,|\xi|)\right)u_0(\cdot)+\ml{F}^{-1}_{\xi\to x}\left(\widehat{K}_1^u(t,|\xi|)\right)u_1(\cdot)+\ml{F}^{-1}_{\xi\to x}\left(\widehat{K}_2^u(t,|\xi|)\right)\theta_0(\cdot)\right]
\end{align*}
in the $\dot{H}^s_q$ sense, by replacing $1-1/r=1/m-1/q$ we summarize all above estimates to complete the desired estimate \eqref{Est-01}. For another, with the aid of  $|\xi|^{-2}\partial_t^2\widehat{K}_j(t,|\xi|)\sim |\xi|^2\widehat{K}_j(t,|\xi|)$ for any $j\in\{0,1,2\}$, we follow the same approach for proving \eqref{Est-01} to derive the estimate for $|D|^s\theta(t,\cdot)$ in the $L^q$ norm immediately.
\end{proof}
\begin{remark}
We actually cannot expect decay estimates for $u(t,\cdot)$ in the $L^q$ norm with the additional $L^m$ regular data if
\begin{align*}
n\leqslant\frac{2mq}{q-m}\ \ \mbox{as well as}\ \ n\in\mb{N}_+.
\end{align*}
For example, in \cite[Theorem 1]{Chen-Ikehata=2023} with the $L^2\cap L^1$ data (i.e. $m=1$, $q=2$), the quantity $\|u(t,\cdot)\|_{L^2}$ satisfies large time optimal growth estimate with the rate $t^{1-\frac{n}{4}}$ when $n=1,2,3$; neither growth nor decay estimate when $n=4$.
\end{remark}
\begin{remark}
The additional $L^q$ regularities for the initial data are to avoid some singularities as $t\to0^+$. Otherwise, by dropping the $L^q$ hypotheses for the initial data, one just can get
\begin{align*}
	\|\,|D|^s\partial_t^k u(t,\cdot)\|_{L^q}&\lesssim t^{-\frac{n}{2}\left(\frac{1}{m}-\frac{1}{q}\right)-\frac{s}{2}-k+1}\left(t^{-1}\|u_0\|_{L^m}+\|u_1\|_{L^m}+\|\theta_0\|_{L^m}\right),\\
	\|\,|D|^s\theta(t,\cdot)\|_{L^q}&\lesssim t^{-\frac{n}{2}\left(\frac{1}{m}-\frac{1}{q}\right)-\frac{s}{2}}\left(t^{-1}\|u_0\|_{L^m}+\|u_1\|_{L^m}+\|\theta_0\|_{L^m}\right),
\end{align*}
with $s\in[0,+\infty)$, $k\in\{0,1\}$ for any $t>0$ and any $1\leqslant m\leqslant q\leqslant +\infty$.
\end{remark}
\begin{remark}\label{Rem-Compare}
Let us compare the known results for the linear thermoelastic plate system \eqref{Linear-TEP} with our result. Remarkably, Theorem \ref{Thm-Linear-u} covers all known estimates in \cite{Munoz-Racke=1995,Denk-Racke=2006,Racke-Ueda=2016,Chen-Ikehata=2023} for $\mb{R}^n$.
\begin{itemize}
	\item The sharp decay rates for the energy term in \cite[Theorem 3.6]{Munoz-Racke=1995} and \cite[Theorem 4.2 with $\eta=2$, $\alpha=\beta=1/2$]{Denk-Racke=2006}, i.e.
	\begin{align*}
		\|(\Delta u,u_t,\theta)(t,\cdot)\|_{(L^q)^3}\lesssim t^{-\frac{n}{2}\left(1-\frac{2}{q}\right)}\|(\Delta u_0,u_1,\theta_0)\|_{(L^{\frac{q}{q-1}})^3}
	\end{align*}
with $q\in[2,+\infty]$, coincide with ours in Theorem \ref{Thm-Linear-u} with $m=q'$.
	\item The sharp decay estimates for the energy term in \cite[Theorem 2.2 and Remark 1]{Racke-Ueda=2016}, i.e.
	\begin{align*}
		\|\,|D|^s(\Delta u,u_t,\theta)(t,\cdot)\|_{(L^2)^3}\lesssim (1+t)^{-\frac{n}{4}-\frac{s}{2}}\|(\Delta u_0,u_1,\theta_0)\|_{(\dot{H}^s\cap L^1)^3}
	\end{align*}
	with $s\geqslant0$, coincide with ours in Theorem \ref{Thm-Linear-u} with $q=2$, $m=1$.
	\item The optimal growth/decay estimates in \cite[Theorem 1]{Chen-Ikehata=2023}, i.e.
	\begin{align*}
		\|u(t,\cdot)\|_{L^2}\simeq t^{1-\frac{n}{4}}\ \ \mbox{and}\ \ \|\theta(t,\cdot)\|_{L^2}\simeq t^{-\frac{n}{4}}
	\end{align*}
	for large time $t\gg1$ provided that $\int_{\mb{R}^n} u_1(x)\,\mathrm{d}x\neq0$ and $\int_{\mb{R}^n} \theta_0(x)\,\mathrm{d}x\neq0$, exactly coincide with ours in Theorem \ref{Thm-Linear-u} with $k=s=0$, $q=2$, $m=1$.
\end{itemize}
For another, concerning the half-space $\mb{R}_+^n:=\{x=(x_1,\dots,x_n):\ x_n>0\}$ the sharp decay rates for the energy term in \cite[Theorem 1.4 and Remark 1.5]{Naito-Shibata=2009}, i.e. 
\begin{align*}
\|\,|D|^s(\Delta u,u_t,\theta)(t,\cdot)\|_{(L^q(\mb{R}^n_+))^3}\lesssim t^{-\frac{n}{2}\left(\frac{1}{m}-\frac{1}{q}\right)-\frac{s}{2}}\|(\Delta u_0,u_1,\theta_0)\|_{(L^m(\mb{R}^n_+))^3}
\end{align*}
with $1<m\leqslant q\leqslant +\infty$ when $s=0,1$ and $1<m\leqslant q<+\infty$ when $s=2$, coincide with ours in Theorem \ref{Thm-Linear-u}. Note that the endpoint cases with $q\in\{1,+\infty\}$ also have been completed for any $s\geqslant 0$ in our result.
\end{remark}

Combining \eqref{Exist-01}, \eqref{Exist-02}, \eqref{Exist-03} and the representations in the Fourier space, we conclude the next $L^q$ well-posedness for our linear Cauchy problem. Furthermore, according to the regularity-gain-type pointwise estimates in Remark \ref{Rem-Pointwise}, the ultra-analytic smoothing effect arises. Let us recall the generalized Gevrey space (see \cite{Rodino=1993}) via the Fourier transform
\begin{align*}
	\Gamma^{\kappa}:=\left\{g\in L^2:\  \exp\left(c\langle\xi\rangle^{\frac{1}{\kappa}}\right)\widehat{g}\in L^2\right\}\ \ \mbox{with}\ \ \kappa\in(0,+\infty),
\end{align*}
where $\langle\xi\rangle^2:=1+|\xi|^2$ denotes the Japanese bracket. Particularly, the space $\Gamma^{\kappa}$ with $\kappa\in(0,1)$ consists of ultra-analytic functions, with $\kappa=1$ consists of analytic functions, and with $\kappa\in(1,+\infty)$ consists of classical Gevrey functions.
\begin{coro}\label{Coro-Linear-GESDS}
Let $s\in[0,+\infty)$. Suppose that the initial data $(u_0,u_1,\theta_0)\in H^{s+2}_q\times H^s_q\times H^s_q$. Then, there is a pair of uniquely determined Sobolev solutions
\begin{align*}
(u,\theta)\in\big(\ml{C}([0,+\infty),H^{s+2}_q)\cap \ml{C}^1([0,+\infty),H^s_q)\big)\times \ml{C}([0,+\infty),H^s_q)
\end{align*}
for any $q\in[1,+\infty]$, to the linear thermoelastic plate system \eqref{Linear-TEP}. Particularly, in the case $q=2$ the solutions further belong to the ultra-analytic space  such that
\begin{align*}
|D|^{s+2} u(t,\cdot),|D|^su_t(t,\cdot),|D|^s\theta(t,\cdot)\in \Gamma^{\frac{1}{2}}\ \ \mbox{for any}\ \ t>0,
\end{align*}
provided that $|D|^{s+2}u_0,|D|^su_1,|D|^s\theta_0\in L^2$ for any $s\in[0,+\infty)$.
\end{coro}

\section{Global in time solutions for the nonlinear systems}\setcounter{equation}{0}\label{Section-GESDS}
\subsection{Philosophy of our proof}\label{Sub-Philos}
\hspace{5mm}For any $T_*>0$ we introduce the evolution space $\ml{X}(T_*)$ of solution $\ml{U}=\ml{U}(t,x)$ in the vector sense that $\ml{U}:=(u^N,\theta^N)^{\mathrm{T}}$ by
\begin{align*}
	\ml{X}(T_*):=\left(\ml{C}([0,T_*],H^2_1\cap H^2_{\infty})\cap\ml{C}^1([0,T_*],L^1\cap L^{\infty}) \right)\times \ml{C}([0,T_*],L^1\cap L^{\infty})
\end{align*}
equipped the time-weighted norm
\begin{align*}
	\|\ml{U}\|_{\ml{X}(T_*)}&:=\sup\limits_{t\in[0,T_*]}\Big(\sum\limits_{s=0,1,2}\left((1+t)^{\frac{s}{2}-1}\|\,|D|^su^N(t,\cdot)\|_{L^1}+(1+t)^{\frac{n}{2}+\frac{s}{2}-1}\|\,|D|^su^N(t,\cdot)\|_{L^{\infty}}\right)\\
	&\qquad\qquad\quad+\|u_t^N(t,\cdot)\|_{L^1}+(1+t)^{\frac{n}{2}}\|u_t^N(t,\cdot)\|_{L^{\infty}}+\|\theta^N(t,\cdot)\|_{L^1}+(1+t)^{\frac{n}{2}}\|\theta^N(t,\cdot)\|_{L^{\infty}}\Big).
\end{align*}
Its time-dependent weighted functions are strongly motivated by some estimates of solutions to the corresponding linearized Cauchy problem \eqref{Linear-TEP}, i.e. Theorem \ref{Thm-Linear-u} so that we may observe the effect of no loss of estimates later.

We next deduce the formal representations of mild solutions $u^N,\theta^N$ to the nonlinear problem \eqref{General-Semilinear-TEP}, whose coupling structure appears not only in the nonlinear part but also strongly in the linear part. Different from the study for coupled systems of semilinear scalar equations (e.g. \cite{Nishihara-Wakasugi=2014,Nishihara-Wakasugi=2015,Djaouti-Reissig=2018,Djaouti=2018} for weakly coupled nonlinearities, and \cite{Ogawa-Takeda=2010,Ogawa-Takeda=2011} for strongly coupled nonlinearities), an additional analysis is needed since the linear coupled system. 

Let us consider another vector unknown $\ml{V}=\ml{V}(t,x)$ such that $\ml{V}:=(u^N,u_t^N,\theta^N)^{\mathrm{T}}$, which satisfies the first order in time coupled system
\begin{align}\label{System-V-nonlinear}
	\begin{cases}
		\ml{V}_t-\ml{A}\ml{V}=F(\ml{V}),&x\in\mb{R}^n,\ t>0,\\
		\ml{V}(0,x)=\ml{V}_0(x),&x\in\mb{R}^n,
	\end{cases}
\end{align}
with the differential operator $\ml{A}$ and the nonlinearity $F(\ml{V})$ such that
\begin{align*}
	\ml{A}:=\left(
	{\begin{array}{*{20}c}
			0 & 1 & 0\\
			-\Delta^2 & 0 & -\Delta\\
			0 & \Delta & \Delta
	\end{array}}
	\right)\ \ \mbox{and}\ \ F(\ml{V}):=\left(
	{\begin{array}{*{20}c}
			0 \\
			f_1^N\\
			f_2^N\\
	\end{array}}
	\right).
\end{align*}
Clearly, the solution $\ml{V}^{\lin}=\ml{V}^{\lin}(t,x)$ to the corresponding linearized model of \eqref{System-V-nonlinear} with the vanishing right-hand side is given by
\begin{align*}
	\left(
	{\begin{array}{*{20}c}
			u \\
			u_t\\
			\theta\\
	\end{array}}
	\right)=\ml{V}^{\lin}=\mathrm{e}^{\ml{A}t}\,\ml{V}_0^{\lin}=\left(
	{\begin{array}{*{20}c}
			K_0^u(t,|D|) & K_1^u(t,|D|) & K_2^u(t,|D|) \\
			\partial_tK_0^u(t,|D|) & \partial_tK_1^u(t,|D|) & \partial_tK_2^u(t,|D|)\\
			K_0^{\theta}(t,|D|) & K_1^{\theta}(t,|D|) & K_2^{\theta}(t,|D|)\\
	\end{array}}
	\right)\left(
	{\begin{array}{*{20}c}
			u_0 \\
			u_1\\
			\theta_0\\
	\end{array}}
	\right)
\end{align*}
thanks to \eqref{Rep-u-01} as well as \eqref{Rep-theta-01}. It suggests the representation of operator semigroup $\mathrm{e}^{\ml{A}t}$. From Duhamel's principle, concerning the first order in time nonlinear coupled system \eqref{System-V-nonlinear}, its mild solution is formally provided by
\begin{align*}
	\ml{V}(t,x)=\mathrm{e}^{\ml{A}t}\,\ml{V}_0(x)+\int_0^t\mathrm{e}^{\ml{A}(t-\tau)}\,F\big(\ml{V}(\tau,x)\big)\,\mathrm{d}\tau,
\end{align*}
to be specific,
\begin{align}\label{Duhamel-nonlinear}
	\int_0^t\mathrm{e}^{\ml{A}(t-\tau)}\,F\big(\ml{V}(\tau,x)\big)\,\mathrm{d}\tau=\int_0^t\left(
	{\begin{array}{*{20}c}
			K_1^u(t-\tau,|D|)f_1^N(\tau,x) + K_2^u(t-\tau,|D|)f_2^N(\tau,x) \\
			\partial_tK_1^u(t-\tau,|D|)f_1^N(\tau,x) + \partial_tK_2^u(t-\tau,|D|)f_2^N(\tau,x)\\
			K_1^{\theta}(t-\tau,|D|)f_1^N(\tau,x) + K_2^{\theta}(t-\tau,|D|)f_2^N(\tau,x)\\
	\end{array}}
	\right)\mathrm{d}\tau.
\end{align}
For this reason, we are able to construct the following nonlinear integral operator for the crucial unknown $\ml{U}$:
\begin{align*}
	\ml{N}:\ \ml{U}\in\ml{X}(T_*)\to\ml{N}[\ml{U}]:=\ml{U}_{\lin}+\ml{U}_{\nlin}.
\end{align*}
Here, $\ml{U}_{\lin}=\ml{U}_{\lin}(t,x)$ denotes the solution to the linearized thermoelastic plate system \eqref{Linear-TEP}, namely, $\ml{U}_{\lin}:=(u,\theta)^{\mathrm{T}}$. From the representations in \eqref{Duhamel-nonlinear} and $\ml{V}(t,x)$, the nonlinear part $\ml{U}_{\nlin}=\ml{U}_{\nlin}(t,x)$ is expressed via
\begin{align*}
	\ml{U}_{\nlin}(t,x)=\left(
	{\begin{array}{*{20}c}
			u_{\nlin}(t,x) \\
			\theta_{\nlin}(t,x)\\
	\end{array}}
	\right):=\int_0^t\left(
	{\begin{array}{*{20}c}
			K_1^u(t-\tau,|D|)f_1^N(\tau,x) + K_2^u(t-\tau,|D|)f_2^N(\tau,x) \\
			K_1^{\theta}(t-\tau,|D|)f_1^N(\tau,x) + K_2^{\theta}(t-\tau,|D|)f_2^N(\tau,x)\\
	\end{array}}
	\right)\mathrm{d}\tau.
\end{align*}
In other words, $\ml{U}_{\lin}+\ml{U}_{\nlin}=(\ml{V}_{(1)},\ml{V}_{(3)})^{\mathrm{T}}$ with the vector $\ml{V}=(\ml{V}_{(1)},\ml{V}_{(2)},\ml{V}_{(3)})^{\mathrm{T}}$ is the formal mild solution to the nonlinear problem \eqref{General-Semilinear-TEP}. We are going  to rigorously justify this statement in an evolution space $\ml{X}(+\infty)$.

To achieve our main goal, i.e. $\ml{U}$ is a fixed point of operator $\ml{N}$ in $\ml{X}(+\infty)$ in the sense of small size by applying the Banach contraction principle, we should show that the following fundamental inequalities:
\begin{align}\label{Crucial-01}
	\|\ml{N}[\ml{U}]\|_{\ml{X}(T_*)}&\lesssim\|(u_0^N,u_1^N,\theta_0^N)\|_{\ml{D}}+\|\ml{U}\|_{\ml{X}(T_*)}^p,\\
	\|\ml{N}[\ml{U}]-\ml{N}[\widetilde{\ml{U}}]\|_{\ml{X}(T_*)}&\lesssim\|\ml{U}-\widetilde{\ml{U}}\|_{\ml{X}(T_*)}\left(\|\ml{U}\|_{\ml{X}(T_*)}^{p-1}+\|\widetilde{\ml{U}}\|_{\ml{X}(T_*)}^{p-1} \right),\label{Crucial-02}
\end{align}
with $p=\min\{p_1,p_2\}>1$, hold for any $\ml{U},\widetilde{\ml{U}}\in \ml{X}(T_*)$, where all unexpressed multiplicative constants are independent of $T_*$. Note that the estimates in Theorem \ref{Thm-Linear-u} show $\ml{U}_{\lin}\in \ml{X}(T_*)$  for any $T_*>0$ (or see Corollary \ref{Coro-Linear-GESDS} with $s=0$ as well as $q\in\{1,+\infty\}$) and the uniformly in time bounded estimate
\begin{align}\label{Crucial-03}
	\|\ml{U}_{\lin}\|_{\ml{X}(T_*)}\lesssim\|(u_0^N,u_1^N,\theta_0^N)\|_{\ml{D}}
\end{align}
with the data space $\ml{D}=(H^2_1\cap H^2_{\infty})\times (L^1\cap L^{\infty})\times (L^1\cap L^{\infty})$.

\subsection{A priori  $L^r$ estimates for the nonlinearities}\label{Subsection-A prior}
\hspace{5mm} Let us first recall the definition of evolution space $\ml{X}(T_*)$. By applying the Riesz-Thorin interpolation between the $L^1$ and $L^{\infty}$ norms for $|D|^su^N(\tau,\cdot)$, $u^N_t(\tau,\cdot)$, $\theta^N(\tau,\cdot)$ with $s\in\{0,\sigma\}$, respectively, we are able to conclude
\begin{align*}
	\|\,|u^N(\tau,\cdot)|^{p_j}\|_{L^r}&=\|u^N(\tau,\cdot)\|_{L^{rp_j}}^{p_j}\lesssim (1+\tau)^{\left(1-\frac{n}{2}\right)p_j+\frac{n}{2r}}\|\ml{U}\|_{\ml{X}(\tau)}^{p_j},\\
	\|\,|\,|\nabla|^{\sigma}u^N(\tau,\cdot)|^{p_j}\|_{L^r}&=\|\,|\nabla|^{\sigma}u^N(\tau,\cdot)\|_{L^{rp_j}}^{p_j}\lesssim (1+\tau)^{\left(1-\frac{n}{2}-\frac{\sigma}{2}\right)p_j+\frac{n}{2r}}\|\ml{U}\|_{\ml{X}(\tau)}^{p_j},\\
	\|\,|u^N_t(\tau,\cdot)|^{p_j}\|_{L^r}&=\|u^N_t(\tau,\cdot)\|_{L^{rp_j}}^{p_j}\lesssim (1+\tau)^{-\frac{n}{2}p_j+\frac{n}{2r}}\|\ml{U}\|_{\ml{X}(\tau)}^{p_j},
\end{align*}
moreover,
\begin{align*}
	\|\,|\theta^N(\tau,\cdot)|^{p_j}\|_{L^r}=\|\theta^N(\tau,\cdot)\|_{L^{rp_j}}^{p_j}\lesssim (1+\tau)^{-\frac{n}{2}p_j+\frac{n}{2r}}\|\ml{U}\|_{\ml{X}(\tau)}^{p_j},
\end{align*}
for $\sigma\in\{1,2\}$ and any $r\in[1,+\infty]$. Note that for the fractional number $\sigma$ we need to apply the fractional Gagliardo-Nirenberg inequality as said in Remark \ref{Rem-1-infty}, which restricts $rp_j\in(1,+\infty)$, i.e. $r\in[1,+\infty)$. Moreover, the last four a priori estimates are independent of the explicit expression of nonlinearity. It may be applied in a general power nonlinear term.

Consequently, the power nonlinearities  of the additivity-type \eqref{nonlinearity} are estimated by
\begin{align*}
	\|f_j^N(\tau,\cdot)\|_{L^r}&\lesssim (1+\tau)^{-\frac{n}{2}p_j+\frac{n}{2r}}\left(c_{1,j}(1+\tau)^{p_j}+c_{2,j}(1+\tau)^{\frac{2-\sigma}{2}p_j}+c_{3,j}+c_{4,j}\right)\|\ml{U}\|_{\ml{X}(\tau)}^{p_j}\\
	&\lesssim\begin{cases}
		(1+\tau)^{\left(1-\frac{n}{2}\right)p_j+\frac{n}{2r}}\|\ml{U}\|_{\ml{X}(\tau)}^{p_j}&\mbox{if}\ \ c_{1,j}\neq0\\
		(1+\tau)^{\left(1-\frac{n}{2}-\frac{\sigma}{2}\right)p_j+\frac{n}{2r}}\|\ml{U}\|_{\ml{X}(\tau)}^{p_j}&\mbox{if}\ \ c_{1,j}=0\ \ \mbox{but}\ \ c_{2,j}\neq0\\
		(1+\tau)^{-\frac{n}{2}p_j+\frac{n}{2r}}\|\ml{U}\|_{\ml{X}(\tau)}^{p_j}&\mbox{if}\ \ c_{1,j}=c_{2,j}=0\ \ \mbox{but}\ \ c_{3,j}\neq0\ \ \mbox{or} \ \ c_{4,j}\neq0\\
		0&\mbox{if}\ \ c_{1,j}=c_{2,j}=c_{3,j}=c_{4,j}=0
	\end{cases}\\
	&=:\ml{B}_{j,r}(1+\tau)\|\ml{U}\|_{\ml{X}(\tau)}^{p_j}
\end{align*}
for $j\in\{1,2\}$. Here, one may be aware of the influence of coefficients $c_{d,j}$ of nonlinearities on global in time behavior of solutions. Remark that the $\tau$-dependent coefficient $\ml{B}_{j,r}(1+\tau)$ is a decreasing function with respect to $r\in[1,+\infty]$, namely,
\begin{align*}
\|f_j^N(\tau,\cdot)\|_{ L^{\infty}\cap L^1}\lesssim \ml{B}_{j,1}(1+\tau)\|\ml{U}\|_{\ml{X}(\tau)}^{p_j}.
\end{align*}

\subsection{Proof of Theorem \ref{Thm-GESDS}}
\hspace{5mm}Let $s\in\{0,1,2\}$ if $k=0$, and $s=0$ if $k=1$ throughout this subsection. One may employ the derived $L^1-L^1$ estimate, i.e. \eqref{Est-01} with $m=q=1$, in $[0,t]$ to find
\begin{align*}
	&(1+t)^{\frac{s}{2}+k-1}\|\,|D|^s\partial_t^ku_{\nlin}(t,\cdot)\|_{L^1}\\
	&= (1+t)^{\frac{s}{2}+k-1}\left\||D|^s\partial_t^k\int_0^t\big(K_1^u(t-\tau,|D|)f_1^N(\tau,\cdot) + K_2^u(t-\tau,|D|)f_2^N(\tau,\cdot)\big)\,\mathrm{d}\tau \right\|_{L^1}\\
	&\lesssim (1+t)^{\frac{s}{2}+k-1}\int_0^t(1+t-\tau)^{-\frac{s}{2}-k+1}\big(\|f_1^N(\tau,\cdot)\|_{L^1}+\|f_2^N(\tau,\cdot)\|_{L^1}\big)\,\mathrm{d}\tau\\
	&\lesssim \left(\int_0^{t/2}\sum\limits_{j=1,2}\ml{B}_{j,1}(1+\tau)\,\mathrm{d}\tau+(1+t)\sum\limits_{j=1,2}\ml{B}_{j,1}(1+t)\right)\|\ml{U}\|^{\min\{p_1,p_2\}}_{\ml{X}(T_*)}
\end{align*}
since the power $-s/2-k+1\geqslant0$. In the last line, the asymptotic relations $1+t-\tau\approx 1+t$ for $\tau\in[0,t/2]$ and $1+\tau\approx 1+t$ for $\tau\in[t/2,t]$ were used. The integrable condition 
\begin{align}\label{O-01}
\int_0^{+\infty}\ml{B}_{j,1}(1+\tau)\,\mathrm{d}\tau<+\infty\ \ \mbox{for any}\ \ j\in\{1,2\}
\end{align}
 holds due to our assumption \eqref{GESDS-Condition}, which leads to
 \begin{align*}
(1+t)\sum\limits_{j=1,2}\ml{B}_{j,1}(1+t)\lesssim 1\ \ \mbox{for any}\ \ t\geqslant0
 \end{align*}
automatically. For this reason, we arrive at
\begin{align*}
	(1+t)^{\frac{s}{2}+k-1}\|\,|D|^s\partial_t^ku_{\nlin}(t,\cdot)\|_{L^1}\lesssim \|\ml{U}\|^{\min\{p_1,p_2\}}_{\ml{X}(T_*)}.
\end{align*}
Concerning the solution in the $L^{\infty}$ norm, one then employs the $(L^{\infty}\cap L^1)-L^{\infty}$ estimate in $[0,t/2]$ and the $L^{\infty}-L^{\infty}$ estimate in $[t/2,t]$ from Theorem \ref{Thm-Linear-u}, precisely,
\begin{align*}
	&(1+t)^{\frac{n}{2}+\frac{s}{2}+k-1}\|\,|D|^s\partial_t^ku_{\nlin}(t,\cdot)\|_{L^{\infty}}\\
	&\lesssim (1+t)^{\frac{n}{2}+\frac{s}{2}+k-1}\int_0^{t/2}(1+t-\tau)^{-\frac{n}{2}-\frac{s}{2}-k+1}\big(\|f_1^N(\tau,\cdot)\|_{L^{\infty}\cap L^1}+\|f_2^N(\tau,\cdot)\|_{L^{\infty}\cap L^1}\big)\,\mathrm{d}\tau\\
	&\quad+(1+t)^{\frac{n}{2}+\frac{s}{2}+k-1}\int_{t/2}^t(1+t-\tau)^{-\frac{s}{2}-k+1}\big(\|f_1^N(\tau,\cdot)\|_{L^{\infty}}+\|f_2^N(\tau,\cdot)\|_{L^{\infty}}\big)\,\mathrm{d}\tau\\
	&\lesssim \left(\int_0^{t/2}\sum\limits_{j=1,2}\ml{B}_{j,1}(1+\tau)\,\mathrm{d}\tau+(1+t)^{1+\frac{n}{2}}\sum\limits_{j=1,2}\ml{B}_{j,\infty}(1+t)\right)\|\ml{U}\|^{\min\{p_1,p_2\}}_{\ml{X}(T_*)}\\
	&\lesssim \|\ml{U}\|^{\min\{p_1,p_2\}}_{\ml{X}(T_*)},
\end{align*}
where we used the next fact:
\begin{align*}
	(1+t)^{\frac{n}{2}}[\ml{B}_{j,\infty}(1+t)]=\ml{B}_{j,1}(1+t)\ \ \mbox{for}\ \ j\in\{ 1,2\}
\end{align*}
and our assumption \eqref{GESDS-Condition} again leading to \eqref{O-01}. By the same way as those for $\partial_tu_{\nlin}(t,\cdot)$ in the $L^1$ and $L^{\infty}$ norms, respectively, in the above calculation (because all estimates for $u_t$ and $\theta$ are the same in Theorem \ref{Thm-Linear-u}) one claims easily
\begin{align*}
	\|\theta_{\nlin}(t,\cdot)\|_{L^1}+(1+t)^{\frac{n}{2}}\|\theta_{\nlin}(t,\cdot)\|_{L^{\infty}}\lesssim \|\ml{U}\|^{\min\{p_1,p_2\}}_{\ml{X}(T_*)}.
\end{align*}
Let us combine the last estimates to show
\begin{align*}
	\|\ml{U}_{\nlin}\|_{\ml{X}(T_*)}\lesssim \|\ml{U}\|^{\min\{p_1,p_2\}}_{\ml{X}(T_*)},
\end{align*}
which implies our desired estimate \eqref{Crucial-01} by using \eqref{Crucial-03} immediately.

To end our proof, we eventually sketch the idea for justifying the Lipschitz condition \eqref{Crucial-02}. Due to the fact that
\begin{align*}
	|\,|v|^{p_j}-|\widetilde{v}|^{p_j}|\lesssim |v-\widetilde{v}|\left(|v|^{p_j-1}+|\widetilde{v}|^{p_j-1}\right)
\end{align*}
for $v=u^N,|\nabla|^{\sigma}u^N,u_t^N,\theta^N$ and their corresponding $\widetilde{v}$, by applying H\"older's inequality and the Riesz-Thorin interpolation theorem again (similarly to those in Subsection \ref{Subsection-A prior}) the difference of nonlinearities satisfies
\begin{align*}
	\|f_j^N(\tau,\cdot)-\widetilde{f}^N_j(\tau,\cdot)\|_{L^r}\lesssim \ml{B}_{j,r}(1+\tau)\,\|\ml{U}-\widetilde{\ml{U}}\|_{\ml{X}(\tau)}\left(\|\ml{U}\|_{\ml{X}(\tau)}^{p_j-1}+\|\widetilde{\ml{U}}\|_{\ml{X}(\tau)}^{p_j-1} \right)
\end{align*}
for $j\in\{1,2\}$ and any $r\in[1,+\infty]$. With the aid of \eqref{GESDS-Condition}, one may repeat the same calculation as those in deriving \eqref{Crucial-01} to complete the desired estimate \eqref{Crucial-02}.

Hence, associated with the Banach fixed point argument our proof of $\ml{U}\in\ml{X}(+\infty)$ is finished. As a byproduct, for $q\in\{1,+\infty\}$ the definition of evolution space $\ml{X}(+\infty)$ implies
\begin{align*}
	\|\,|D|^s\partial_t^k u^N(t,\cdot)\|_{L^q}&\lesssim(1+t)^{-\frac{n}{2}\left(1-\frac{1}{q}\right)-\frac{s}{2}-k+1}\|(u_0^N,u_1^N,\theta_0^N)\|_{\ml{D}},\\
	\| \theta^N(t,\cdot)\|_{L^q}&\lesssim(1+t)^{-\frac{n}{2}\left(1-\frac{1}{q}\right)}\|(u_0^N,u_1^N,\theta_0^N)\|_{\ml{D}},
\end{align*}
hold for any $t\geqslant0$ with $s\in\{0,1,2\}$ if $k=0$, and $s=0$ if $k=1$. One finally may apply the Riesz-Thorin interpolation between the $L^1$ and $L^{\infty}$ spaces to complete the regularity of solutions in \eqref{Solution-space-Lq} and our desired $L^q$ estimates in Theorem \ref{Thm-GESDS} for any $q\in[1,+\infty]$.

\begin{remark}
The proof of Corollary \ref{Coro-GESDS=2} is similar to the above, where we just need to change some estimates to the new nonlinearities \eqref{strong-nonlinearity}. Let us consider the same evolution space $\ml{X}(T_*)$ as the one in Subsection \ref{Sub-Philos}, and recall the a priori $L^r$ estimates in Subsection \ref{Subsection-A prior}. Then, an application of generalized H\"older's inequality for any $r\in[1,+\infty]$ shows that the nonlinearity $f_1^N$ defined in \eqref{strong-nonlinearity} is estimated by
\begin{align*}
	\|f_1^N(\tau,\cdot)\|_{L^r}&=\||u^N(\tau,\cdot)|^{p_1}\,|\,|D|^{\sigma}u^N(\tau,\cdot)|^{p_2}\,|u^N_t(\tau,\cdot)|^{p_3}\,|\theta^N(\tau,\cdot)|^{p_4} \|_{L^r}\\
	&\leqslant\|\,|u^N(\tau,\cdot)|^{p_1}\|_{L^{r_1}}\|\,|\,|D|^{\sigma}u^N(\tau,\cdot)|^{p_2}\|_{L^{r_2}}\|\,|u^N_t(\tau,\cdot)|^{p_3}\|_{L^{r_3}}\|\,|\theta^N(\tau,\cdot)|^{p_4}\|_{L^{r_4}}\\
	&\lesssim (1+\tau)^{\left(1-\frac{n}{2}\right)p_1+\frac{n}{2r_1}+\left(1-\frac{n}{2}-\frac{\sigma}{2}\right)p_2+\frac{n}{2r_2}-\frac{n}{2}p_3+\frac{n}{2r_3}-\frac{n}{2}p_4+\frac{n}{2r_4}}\|\ml{U}\|_{\ml{X}(\tau)}^{p_1+p_2+p_3+p_4}\\
	&\lesssim (1+\tau)^{p_1+\left(1-\frac{\sigma}{2}\right)p_2-\frac{n}{2}\sum\limits_{\ell=1,\dots,4}p_\ell+\frac{n}{2r}}\|\ml{U}\|_{\ml{X}(\tau)}^{p_1+p_2+p_3+p_4}
\end{align*}
with $1/r_1+1/r_2+1/r_3+1/r_4=1/r$ and $r_{\ell}\in[1,+\infty]$. Analogously,
\begin{align*}
	\|f_2^N(\tau,\cdot)\|_{L^r}\lesssim (1+\tau)^{\bar{p}_1+\left(1-\frac{\sigma}{2}\right)\bar{p}_2-\frac{n}{2}\sum\limits_{\ell=1,\dots,4}\bar{p}_\ell+\frac{n}{2r}}\|\ml{U}\|_{\ml{X}(\tau)}^{\bar{p}_1+\bar{p}_2+\bar{p}_3+\bar{p}_4}
\end{align*}
for any $r\in[1,+\infty]$. By our assumptions on the powers, the following estimates are valid:
\begin{align*}
\int_0^{t/2}\|f_1^N(\tau,\cdot)\|_{L^1}\,\mathrm{d}\tau&\lesssim \int_0^{t/2}(1+\tau)^{p_1+\left(1-\frac{\sigma}{2}\right)p_2-\frac{n}{2}\sum\limits_{\ell=1,\dots,4}p_\ell+\frac{n}{2}}\,\mathrm{d}\tau\,\|\ml{U}\|_{\ml{X}(T_*)}^{p_1+p_2+p_3+p_4}\lesssim \|\ml{U}\|_{\ml{X}(T_*)}^{p_1+p_2+p_3+p_4},\\
\int_0^{t/2}\|f_2^N(\tau,\cdot)\|_{L^1}\,\mathrm{d}\tau&\lesssim \int_0^{t/2}(1+\tau)^{\bar{p}_1+\left(1-\frac{\sigma}{2}\right)\bar{p}_2-\frac{n}{2}\sum\limits_{\ell=1,\dots,4}\bar{p}_\ell+\frac{n}{2}}\,\mathrm{d}\tau\,\|\ml{U}\|_{\ml{X}(T_*)}^{\bar{p}_1+\bar{p}_2+\bar{p}_3+\bar{p}_4}\lesssim \|\ml{U}\|_{\ml{X}(T_*)}^{\bar{p}_1+\bar{p}_2+\bar{p}_3+\bar{p}_4},
\end{align*}
instead of \eqref{O-01}. Consequently, the estimates \eqref{Crucial-01} and \eqref{Crucial-02} hold with
\begin{align*}
p=\min\{p_1+p_2+p_3+p_4,\bar{p}_1+\bar{p}_2+\bar{p}_3+\bar{p}_4\}>1.
\end{align*}  It immediately concludes the global in time existence result in Corollary \ref{Coro-GESDS=2}.
\end{remark}

\section{Blow-up phenomena for the nonlinear systems}\setcounter{equation}{0}\label{Section-Blow-up}
\hspace{5mm}Let us define the test function $\eta(t)\in\ml{C}_0^{\infty}([0,+\infty))$ and the radially symmetric test function $\varphi(x)\in\ml{C}_0^{\infty}(\mb{R}^n)$ such that
\begin{align*}
\eta(t):=\begin{cases}
1&\mbox{if}\ \ t\in[0,1/2]\\
\searrow&\mbox{if}\ \ t\in(1/2,1)\\
0&\mbox{if}\ \ t\in[1,+\infty)
\end{cases}\ \ \mbox{and}\ \ 
\varphi(|x|):=\begin{cases}
1&\mbox{if}\ \ |x|\in[0,1/2]\\
\searrow&\mbox{if}\ \ |x|\in(1/2,1)\\
0&\mbox{if}\ \ |x|\in[1,+\infty)
\end{cases}
\end{align*}
according to the radial function $\varphi(x)\equiv\varphi(|x|)$. Here, the notation $\searrow$ means the non-increasing behavior.
Additionally, they are assumed to satisfy the following bounded estimates (see, for example, \cite{Mitidieri-Poho=2001,D'Ambro-Lucen=2003}):
\begin{align}\label{test function eta varphi 1}
		[\eta(t)]^{-\frac{r'}{r}}\Big(|\eta'(t)|^{r'}+|\eta''(t)|^{r'}\Big)&\leqslant C,\qquad[\varphi(x)]^{-\frac{r'}{r}}\Big(|\Delta\varphi(x)|^{r'}+|\Delta^2\varphi(x)|^{r'}\Big)\leqslant C,
\end{align}
with the exponent $r\in\{p_1,p_2\}$ and the H\"older's conjugate $r'=r/(r-1)$. These test functions are strongly motivated by the pioneering work \cite{Zhang=2001} for the semilinear classical damped wave equation, and its applications in \cite{Wakasugi=2014,D'Abbicco-Reissig=2014,Pham-Kainae-Reissig=2015,do-Pa-Rei=2017,D'Abbicco-Ebert=2017} for several semilinear evolution models with different nonlinearities.

\subsection{Proof of Theorem \ref{Thm-Blow-up-01}}
\hspace{5mm}For any $R>0$ one may introduce the following test function with separate variables:
\begin{align*}
	\psi_R(t,x):=\eta_R(t)\,\varphi_R(x):=\eta(R^{-2}t)\,\varphi(R^{-1}x).
\end{align*}
Let $\Phi_R(t)\in\ml{C}_0^{\infty}([0,+\infty))$ be another test function such that
\begin{align*}
	\Phi_R(t):=\int_t^{+\infty}\eta_R(\tau)\,\mathrm{d}\tau,
\end{align*}
which fulfills $\Phi'_R(t)=-\eta_R(t)$ and $\Phi_R(t)\leqslant \Phi_R(0)$ for any $t\geqslant0$.

Choosing the test functions by $\Psi_1(t,x)=\Psi_2(t,x)=\psi_R(t,x)$ in Definition \ref{Defn-weak-01}, with $u_0^N(x)=0$ one knows
	\begin{align}\label{est 3.1}
	&\int_0^{+\infty}\int_{\mb{R}^n}f_1^N(t,x)\,\psi_R(t,x)\,\mathrm{d}x\,\mathrm{d}t+\int_{\mb{R}^n}u^N_1(x)\,\varphi_R(x)\,\mathrm{d}x\notag\\
	&=\int_0^{+\infty}\int_{\mb{R}^n}u^N_t(t,x)\big(-\eta_R'(t)\,\varphi_R(x)+\Phi_R(t)\Delta^2\varphi_R(x)\big)\,\mathrm{d}x\,\mathrm{d}t\notag\\
	&\quad+\int_0^{+\infty}\int_{\mb{R}^n}\theta^N(t,x)\,\eta_R(t)\Delta \varphi_R(x)\,\mathrm{d}x\,\mathrm{d}t
\end{align}
as well as
\begin{align}\label{est 3.2}
	&\int_0^{+\infty}\int_{\mb{R}^n}f_2^N(t,x)\,\psi_R(t,x)\,\mathrm{d}x\,\mathrm{d}t+\int_{\mb{R}^n}\theta_0^N(x)\,\varphi_R(x)\,\mathrm{d}x\notag\\
	&=\int_0^{+\infty}\int_{\mb{R}^n}\theta^N(t,x)\big(-\eta_R'(t)\,\varphi_R(x)-\eta_R(t)\Delta \varphi_R(x)\big)\,\mathrm{d}x\,\mathrm{d}t\notag\\
	&\quad-\int_0^{+\infty}\int_{\mb{R}^n}u^N_t(t,x)\, \eta_R(t)\Delta\varphi_R(x)\,\mathrm{d}x\,\mathrm{d}t.
\end{align}
Thanks to our assumptions $c_{3,1}+c_{3,2}>0$ and $c_{4,1}+c_{4,2}>0$, the nonlinear terms can be controlled from the below by
\begin{align*}
	&\int_0^{+\infty}\int_{\mb{R}^n}\big(f_1^N(t,x)+f_2^N(t,x)\big)\psi_R(t,x)\,\mathrm{d}x\,\mathrm{d}t\\
	&\geqslant \sum\limits_{j=1,2}\int_0^{+\infty}\int_{\mb{R}^n}\big(c_{3,j}|u^N_t(t,x)|^{p_j}\psi_R(t,x)+c_{4,j}|\theta^N(t,x)|^{p_j}\psi_R(t,x)\big)\,\mathrm{d}x\,\mathrm{d}t=:\sum\limits_{j=1,2}(c_{3,j}I_{R,j}+c_{4,j}J_{R,j})
\end{align*}
since $c_{d,j}\geqslant0$ for any $d\in\{1,\dots,4\}$ and $j\in\{1,2\}$. Furthermore, there exist $\ell_1,\ell_2\in\{1,2\}$ such that $c_{3,\ell_1}>0$ and $c_{4,\ell_2}>0$. Consequently,
\begin{align*}
	\sum\limits_{j=1,2}(c_{3,j}I_{R,j}+c_{4,j}J_{R,j})\geqslant c_{3,\ell_1}I_{R,\ell_1}+c_{4,\ell_2}J_{R,\ell_2}.
\end{align*}
The combination of \eqref{est 3.1} and \eqref{est 3.2} yields
\begin{align*}
	&c_{3,\ell_1}I_{R,\ell_1}+c_{4,\ell_2}J_{R,\ell_2}+\int_{\mb{R}^n}\big(u_1^N(x)+\theta_0^N(x)\big)\varphi_R(x)\,\mathrm{d}x\\
	&\leqslant \int_0^{+\infty}\int_{\mb{R}^n}u^N_t(t,x)\big(-\eta'_R(t)\,\varphi_R(x)-\eta_R(t)\Delta\varphi_R(x)+\Phi_R(t)\Delta^2\varphi_R(x)\big)\,\mathrm{d}x\,\mathrm{d}t\\
	&\quad-\int_0^{+\infty}\int_{\mb{R}^n}\theta^N(t,x)\,\eta'_R(t)\,\varphi_R(x)\,\mathrm{d}x\,\mathrm{d}t.
\end{align*}
Then, we apply Young's inequality to get
\begin{align*}
	&\frac{c_{3,\ell_1}}{p_{\ell_1}'}I_{R,\ell_1}+\frac{c_{4,\ell_2}}{p_{\ell_2}'}J_{R,\ell_2}+\int_{\mb{R}^n}\big(u_1^N(x)+\theta_0^N(x)\big)\varphi_R(x)\,\mathrm{d}x\\
	&\lesssim\int_0^{+\infty}\int_{\mb{R}^n}\big(\eta_R(t)\,\varphi_R(x)\big)^{-\frac{p'_{\ell_1}}{p_{\ell_1}}}\big(|\eta'_R(t)\,\varphi_R(x)|+|\eta_R(t)\Delta\varphi_R(x)|+|\Phi_R(t)\Delta^2\varphi_R(x)|\big)^{p'_{\ell_1}}\,\mathrm{d}x\,\mathrm{d}t\\
	&\quad+\int_0^{+\infty}\int_{\mb{R}^n}\big(\eta_R(t)\,\varphi_R(x)\big)^{-\frac{p'_{\ell_2}}{p_{\ell_2}}}|\eta'_R(t)\,\varphi_R(x)|^{p'_{\ell_2}}\,\mathrm{d}x\,\mathrm{d}t.
\end{align*}
According to the facts that
\begin{align*}
	\eta_R'(t)&=R^{-2}\,\eta'(R^{-2}t),\ \ \ \  \Delta\varphi_R(x)=R^{-2}\Delta\varphi(R^{-1}x),\\
	\Delta^2\varphi_R(x)&=R^{-4}\Delta^2\varphi(R^{-1}x),\quad\, \Phi_R(t)\leqslant\Phi_R(0)\lesssim R^2,
\end{align*}
we may derive
\begin{align}
		\int_0^{+\infty}\int_{\mb{R}^n}\big(\eta_R(t)\,\varphi_R(x)\big)^{-\frac{r'}{r}}|\eta'_R(t)\,\varphi_R(x)|^{r'}\,\mathrm{d}x\,\mathrm{d}t&\lesssim R^{-2r'+n+2},\label{Est-06}\\
\int_0^{+\infty}\int_{\mb{R}^n}\big(\eta_R(t)\,\varphi_R(x)\big)^{-\frac{r'}{r}}|\eta_R(t)\Delta\varphi_R(x)|^{r'}\,\mathrm{d}x\,\mathrm{d}t&\lesssim R^{-2r'+n+2},\notag\\
\int_0^{+\infty}\int_{\mb{R}^n}\big(\eta_R(t)\,\varphi_R(x)\big)^{-\frac{r'}{r}}|\Phi_R(t)\Delta^2\varphi_R(x)|^{r'}\,\mathrm{d}x\,\mathrm{d}t&\lesssim R^{-2r'+n+2},\notag
\end{align}
for $r\in\{p_{\ell_1},p_{\ell_2}\}$, in which we employed \eqref{test function eta varphi 1} and the supports for $\eta(t),\varphi(x)$. In other words, the crucial estimate occurs
\begin{align}\label{Est-05}
\frac{c_{3,\ell_1}}{p_{\ell_1}'}I_{R,\ell_1}+\frac{c_{4,\ell_2}}{p_{\ell_2}'}J_{R,\ell_2}+\int_{\mb{R}^n}\big(u_1^N(x)+\theta_0^N(x)\big)\varphi_R(x)\,\mathrm{d}x\lesssim R^{-2p'_{\ell_1}+n+2}+R^{-2p'_{\ell_2}+n+2},
\end{align}
in which we considered the linear combination of $I_{R,\ell_1}$ and $J_{R,\ell_2}$ as a whole functional due to its linear coupled system on the left-hand sides of \eqref{General-Semilinear-TEP}.

Let us assume that the pair of solutions $(u^N,\theta^N)$ and, in turn, the functionals $I_{R,\ell_1},J_{R,\ell_2}$ for the semilinear Cauchy problem \eqref{General-Semilinear-TEP} are globally in time defined. Due to our assumption on the initial data, letting $R\to+\infty$ we may derive $\frac{c_{3,\ell_1}}{p_{\ell_1}'}I_{R,\ell_1}+\frac{c_{4,\ell_2}}{p_{\ell_2}'}J_{R,\ell_2}<0$ when 
\begin{align*}
-2p_{\ell_1}'+n+2<0\ \ \mbox{as well as}\ \ -2p_{\ell_2}'+n+2<0,
\end{align*}
simultaneously, that is to say,
\begin{align}\label{A-condition}
	1<p_{\ell_1},p_{\ell_2}<1+\frac{2}{n}\ \ \mbox{for}\  \ \ell_1,\ell_2\in\{1,2\}.
\end{align}
This produces a contradiction via $I_{R,\ell_1},J_{R,\ell_2}\geqslant0$. For this reason, under the condition
\begin{align}\label{B-condition}
1<p_j<1+\frac{2}{n}\ \ \mbox{if}\ \ c_{3,j}>0\ \ \mbox{or}\ \ c_{4,j}>0,
\end{align}
 the global in time weak solution (according to Definition \ref{Defn-weak-01}) does not exist. It is worth mentioning that \eqref{A-condition} is equivalent to \eqref{B-condition} because $c_{3,j}>0\Leftrightarrow \ell_1=j$ and  $c_{4,j}>0\Leftrightarrow \ell_2=j$ for $j\in\{1,2\}$. The critical case
 \begin{align*}
	p_{\ell_1},p_{\ell_2}=1+\frac{2}{n}\ \ \mbox{for}\  \ \ell_1,\ell_2\in\{1,2\}
\end{align*}
can be treated in the standard way (see \cite{Zhang=2001,Wakasugi=2014,D'Abbicco-Reissig=2014} and references therein) associated with the derived estimate \eqref{Est-05}, but we omit the details for the sake of brevity. Therefore, our proof of Theorem \ref{Thm-Blow-up-01} is completed.

\subsection{Proof of Theorem \ref{Thm-Blow-up-02}}
\hspace{5mm}Again, taking our test functions by $\Psi_1(t,x)=\Psi_2(t,x)=\psi_R(t,x)$ in Definition \ref{Defn-weak-02}, with $u_0^N(x)=0$ one may immediately derive
\begin{align}\label{est 3.3}
	&\int_0^{+\infty}\int_{\mb{R}^n}f_1^N(t,x)\,\psi_R(t,x)\,\mathrm{d}x\,\mathrm{d}t+\int_{\mb{R}^n}u_1^N(x)\,\varphi_R(x)\,\mathrm{d}x\notag\\
	&=\int_0^{+\infty}\int_{\mb{R}^n}u^N(t,x)\big(\eta_R''(t)\,\varphi_R(x)+\eta_R(t)\Delta^2\varphi_R(x)\big)\,\mathrm{d}x\,\mathrm{d}t\notag\\
	&\quad+\int_0^{+\infty}\int_{\mb{R}^n}\theta^N(t,x)\,\eta_R(t)\Delta \varphi_R(x)\,\mathrm{d}x\,\mathrm{d}t
\end{align}
as well as
\begin{align}\label{est 3.4}
	&\int_0^{+\infty}\int_{\mb{R}^n}f_2^N(t,x)\,\psi_R(t,x)\,\mathrm{d}x\,\mathrm{d}t+\int_{\mb{R}^n}\theta_0^N(x)\,\varphi_R(x)\,\mathrm{d}x\notag\\
	&=\int_0^{+\infty}\int_{\mb{R}^n}\theta^N(t,x)\big(-\eta'_R(t)\,\varphi_R(x)-\eta_R(t)\Delta \varphi_R(x)\big)\,\mathrm{d}x\,\mathrm{d}t\notag\\
	&\quad+\int_0^{+\infty}\int_{\mb{R}^n}u^N(t,x)\,\eta'_R(t)\Delta\varphi_R(x)\,\mathrm{d}x\,\mathrm{d}t.
\end{align}
Under our conditions $c_{1,1}+c_{1,2}>0$ and $c_{4,1}+c_{4,2}>0$, there exist $\ell_3,\ell_4\in\{1,2\}$ such that $c_{1,\ell_3}>0$ and $c_{4,\ell_4}>0$. Hence, it holds that
\begin{align*}
	&\int_0^{+\infty}\int_{\mb{R}^n}\big(f_1^N(t,x)+f_2^N(t,x)\big)\psi_R(t,x)\,\mathrm{d}x\,\mathrm{d}t\\
	&\geqslant \sum\limits_{j=1,2}\int_0^{+\infty}\int_{\mb{R}^n}\big(c_{1,j}|u^N(t,x)|^{p_j}\psi_R(t,x)+c_{4,j}|\theta^N(t,x)|^{p_j}\psi_R(t,x)\big)\,\mathrm{d}x\,\mathrm{d}t\\
	&=:\sum\limits_{j=1,2}(c_{1,j}L_{R,j}+c_{4,j}J_{R,j})\geqslant c_{1,\ell_3}L_{R,\ell_3}+c_{4,\ell_4}J_{R,\ell_4}.
\end{align*}
Using the last lower bound estimate, the direct sum of \eqref{est 3.3} and \eqref{est 3.4} leads to
\begin{align}
&\frac{c_{1,\ell_3}}{p'_{\ell_3}}L_{R,\ell_3}+\frac{c_{4,\ell_4}}{p'_{\ell_4}}J_{R,\ell_4}+\int_{\mb{R}^n}\big(u_1^N(x)+\theta_0^N(x)\big)\,\varphi_R(x)\,\mathrm{d}x\notag\\
&\lesssim\int_0^{+\infty}\int_{\mb{R}^n}\big(\eta_R(t)\,\varphi_R(x)\big)^{-\frac{p'_{\ell_3}}{p_{\ell_3}}}\big(|\eta''_R(t)\,\varphi_R(x)|+|\eta'_R(t)\Delta\varphi_R(x)|+|\eta_R(t)\Delta^2\varphi_R(x)|\big)^{p'_{\ell_3}}\,\mathrm{d}x\,\mathrm{d}t\notag\\
&\quad\ +\int_0^{+\infty}\int_{\mb{R}^n}\big(\eta_R(t)\,\varphi_R(x)\big)^{-\frac{p'_{\ell_4}}{p_{\ell_4}}}|\eta'_R(t)\,\varphi_R(x)|^{p'_{\ell_4}}\,\mathrm{d}x\,\mathrm{d}t\notag\\
&\lesssim R^{-4p'_{\ell_3}+n+2}+R^{-2p'_{\ell_4}+n+2},\label{Est-07}
\end{align}
where we employed Young's inequality, \eqref{Est-06} and
\begin{align*}
\int_0^{+\infty}\int_{\mb{R}^n}\big(\eta_R(t)\,\varphi_R(x)\big)^{-\frac{r'}{r}}\big|[\mathrm{d}_t^{m_1}\eta_R(t)]\,[\Delta^{2-m_1}\varphi_R(x)]\big|^{r'}\,\mathrm{d}x\,\mathrm{d}t\lesssim R^{-4r'+n+2}
\end{align*}
for $m_1\in\{0,1,2\}$ with the aid of \eqref{test function eta varphi 1}.

Lastly, following the ways of proving Theorem \ref{Thm-Blow-up-01} we may obtain contradictions in the crucial estimate \eqref{Est-07} when 
\begin{align*}
	1<p_{\ell_3}\leqslant1+\frac{4}{(n-2)_+}\ \ \mbox{and}\ \ 1<p_{\ell_4}\leqslant1+\frac{2}{n}\ \ \mbox{for}\  \ \ell_3,\ell_4\in\{1,2\}.
\end{align*}
Thanks to our assumption \eqref{Blow-up-Condition-02}, the proof is completed.

%

\section*{Acknowledgments} 
 Halit S. Aslan is supported in part by Grant 2021/01743-3 from S\~ao Paulo Research Foundation (FAPESP). Wenhui Chen is supported in part by the National Natural Science Foundation of China (grant No. 12301270, grant No. 12171317), 2024 Basic and Applied Basic Research Topic--Young Doctor Set Sail Project (grant No. 2024A04J0016), Guangdong Basic and Applied Basic Research Foundation (grant No. 2023A1515012044). The authors thank Michael Reissig (TU Bergakademie Freiberg) for his suggestions on the global existence part when the authors stayed in Freiberg during 2017-2020.

\end{document}